\documentstyle[amsfonts,latexsym,amsmath,amsthm,amssymb,11pt,leqno]{article}
\begin{document}

\title{\bf Two-parameter Quantum Groups of Exceptional Type $E$-Series and
Convex PBW-Type Basis \footnotetext{$^\star$Corresponding author,
nhhu@math.ecnu.edu.cn}
\thanks{Supported in part by the NNSF (Grant Nos. 10201015, 10431040), the TRAPOYT and the
FUDP from the MOE of China, the SRSTP from the STCSM. }}
 \author{\bf Xiaotang Bai$^{1,2}$   \quad and \quad  Naihong Hu$^{1,\star}$
 \\ {\small $^1$Department of Mathematics, East China Normal University,}\\
 {\small Shanghai 200062, PR. China}\\
 {\small $^2$Department of Mathematics, Harbin University of Science and Technology,}\\
 {\small Harbin 150080, PR. China}}
\date{}
\maketitle
\begin{abstract} The presentation of two-parameter quantum groups of type $E$-series
in the sense of Benkart-Witherspoon [BW1] is given, which has a
Drinfel'd quantum double structure. The universal $R$-matrix and a
convex PBW-type basis are described for type $E_6$ (as a sample),
and the conditions of those isomorphisms from these quantum groups
into the one-parameter quantum doubles are discussed.
\smallskip


 \end{abstract}

\bigskip
\centerline{\textsc{Introduction}}

\bigskip
Two-parameter or multiparameter quantum groups have been
investigated by many authors (see the references in [BW1], [BGH1],
etc.). From another viewpoint based on the work on down-up algebras
(see [B]), Benkart and Witherspoon [BW1] recovered the structure of
two-parameter quantum enveloping algebras of the general linear Lie
algebra $\mathfrak{gl}_n$ and the special linear Lie algebra
$\mathfrak{sl}_n$, which was earlier gotten by Takeuchi [T]. They
studied their finite-dimensional weight representation theory in the
case when $rs^{-1}$ is not a root of unity ([BW2]) and the
restricted quantum version (or say, the small quantum groups in the
two-parameter setting) at $rs^{-1}$ being a root of unity ([BW3]).
Inspired by their work, the two-parameter quantum groups in the
sense of Benkart-Witherspoon corresponding to the orthogonal Lie
algebras $\mathfrak{so}_{2n+1}$ or $\mathfrak{so}_{2n}$ and the
symplectic Lie algebras $\mathfrak{sp}_{2n}$, as well as the
exceptional type $G_2$ were further obtained by Bergeron-Gao-Hu
[BGH1] and Hu-Shi [HS], respectively. Their finite-dimensional
weight representation theory and Lusztig symmetries' property were
systematically established in [BGH2] and [HS]. Actually, this kind
of Lusztig symmetries' property existing from these quantum groups
to their associated objects also reveals the difference with the
standard Drinfeld-Jimbo quantum groups in the one-parameter setting
(see [Ja]).

The aim of this paper is to give the presentation of two-parameter
quantum groups of exceptional type $E$-series, to describe the
universal $R$-matrix and a convex PBW-type basis in terms of Lyndon
words (cf. [LO]), as well as to study those isomorphisms' conditions
from these quantum groups into the one-parameter quantum doubles.
Here we will give a general formalism (see Section 1) of the
presentation of their structural constants, which is actually
applied to all simply-laced types (including types $A$, $D$).

Let $\mathfrak{g}$ denote one of Lie algebras of type $E_6$, $E_7$,
or $E_8$, and $U_{r,s}(\mathfrak{g})$, the two-parameter quantum
enveloping algebra of $\mathfrak{g}$. For simplicity, we will only
write down the results of $E_6$ in this paper (and those for $E_7$ and $E_8$ are similar to be obtained). \\

\centerline{\textsc{1. Presentation of two-parameter quantum group
of type $E$}}
\bigskip
Consider the root system of $E_6$ as a root subsystem of $E_8$.
Assume $\Phi$ is a finite root system of type $E_6$ with a base of
simple roots $\Pi$. We regard $\Phi$ as a subset of a Euclidean
space ${\mathbb{R}}^8$ with an inner product $(\, , \,)$. Let
$\epsilon_1,\epsilon_2,\cdots,\epsilon_8$ denote an orthonormal
basis of ${\mathbb{R}}^8$, and suppose
$\Pi=\{\alpha_1=\frac{1}{2}(\epsilon_1+\epsilon_8)-
\frac{1}{2}(\epsilon_2+\cdots+\epsilon_7),\,
\alpha_2=\epsilon_1+\epsilon_2,\,
\alpha_j=\epsilon_{j-1}-\epsilon_{j-2}\mid 3\leq j\leq 6\}$ and
 $\Phi=\{ \pm(\epsilon_i\pm \epsilon_j)\mid 1\leq j \neq i \leq 5\}\cup
\{\pm \frac{1}{2}(\epsilon_8-\epsilon_7-\epsilon_6+\sum_{i=1}^5\pm
\epsilon_i)\mid \textrm{even number of minus signs }\}$.

Fix two nonzero elements $r,\, s$ in a field $\mathbb{K}$ with
$r\neq s$.

\smallskip
Let $U=U_{r,s}(\mathfrak{g})$ be the unital associative algebra over
$\mathbb{K}$ generated by elements $e_j,f_j,\omega_{i}^{\pm
1},{\omega'}_{i}^{\pm 1} (1\leq i\leq 6)$, which satisfy the
following relations:

\smallskip
(E1) \ $[\omega_{i}^{\pm 1},\omega_{j}^{\pm
1}]=0=[{\omega_{i}'}^{\pm 1}, {\omega_{j}'}^{\pm
1}]=[\omega_{i}^{\pm 1},{\omega_{j}'}^{\pm 1}]$, $\
\omega_{i}\omega_{i}^{-1}=\omega_{j}' {\omega_{j}}'^{-1}=1$.

\smallskip
(E2) \ For $1\leq i,j\leq 6$, we have
$$\omega_{i}e_j\omega_{i}^{-1}=r^{p_{ij}} (s^{-1})^{q_{ij}}e_j, \qquad
\omega_{i}f_j\omega_{i}^{-1}=(r^{-1})^{p_{ij}} s^{q_{ij}}f_j,$$
where $p_{ij}+q_{ij}=(\alpha_i,\alpha_j)$, $p_{ij}, q_{ij} \in \{
0,\pm1\}$, and if $(\alpha_i,\alpha_j)\ne0$, then $p_{ij}-q_{ij}$,
$j-i$ have the same sign.

\smallskip
(E3) \ ${\omega'}_{i}e_j {\omega'}_{i}^{-1}=s^{p_{ij}}
(r^{-1})^{q_{ij}}e_j$, \qquad
${\omega'}_{i}f_j{\omega'}_{i}^{-1}=(s^{-1})^{p_{ij}}
r^{q_{ij}}f_j$.

\smallskip
(E4) \ For $1\leq i,j\leq 6$, we have
$$[e_i,f_j]=\frac{\delta_{i,j}}{r-s}(\omega_i-\omega'_i).$$

(E5) \ For $1\leq i,j\leq 6$, and $(\alpha_i,\alpha_j)=0$,
 $$[e_i,e_j]=[f_i,f_j]=0.$$

(E6) \ For $1\leq i<j\leq 6$, with $a_{ij}=-1$, we have
\begin{gather*}
e_i^2e_j-(r+s)e_ie_je_i+(rs)e_je_i^2=0,\\
e_j^2e_i-(r^{-1}+s^{-1})e_je_ie_j+(r^{-1}s^{-1})e_ie_j^2=0.
\end{gather*}

(E7) \ For $1\leq i<j\leq 6$, with $a_{ij}=-1$, we have
\begin{gather*}
f_jf_i^2-(r+s)f_if_jf_i+(rs)f_i^2f_j=0,\\
f_if_j^2-(r^{-1}+s^{-1})f_jf_if_j+(r^{-1}s^{-1})f_j^2f_i=0.
\end{gather*}

\noindent\textbf{Remark.} \ It is easy to see that when
$(\alpha_i,\alpha_j)=0$, we have two solutions of the equation
$p_{ij}+q_{ij}=(\alpha_i,\alpha_j)$, that is, $p_{ij}=q_{ij}=0$ and
$p_{ij}=\pm 1, q_{ij}=\mp 1$. We have checked that both of them
work, but later on in the next section we only discuss the case when
$p_{ij}=
 q_{ij}=0$ for simplicity. Then for any fixed $(i,j)$,
 $p_{ij}$ and $q_{ij}$ can be determined uniquely.

\medskip
\noindent\textbf{Lemma 1.1.} \textit{For any simply-laced simple Lie
algebra, there hold identities$:$} $p _{ij}=q_{ji}.$

\medskip
\noindent\textit{Proof.} Notice that
$p_{ij}+q_{ij}=(\alpha_i,\alpha_j), \quad
p_{ji}+q_{ji}=(\alpha_j,\alpha_i)$. Since
$(\alpha_i,\alpha_j)=(\alpha_j,\alpha_i)$, then $\{p_{ij},q_{ij}\}$
and $\{p_{ji},q_{ji}\}$ are all the solution of the same equation.
Assume that $i>j$ are two fixed integers, then $\{p_{ij}\leq
q_{ij}\}$ and $\{p_{ji} \geq q_{ji}\}$. Since the solution is
determined uniquely, we can deduce that
$\{p_{ij},q_{ij}\}=\{p_{ji},q_{ji}\}$ and $p _{ij}=q_{ji},q
_{ij}=p_{ji}$. So we get the result. $\qquad \Box$\\

Let $\mathcal{B}=B(\mathfrak{g})$ (resp.
$\mathcal{B}'=B'(\mathfrak{g})$) denote the Hopf subalgebra of
$U=U_{r,s}(\mathfrak{g})$, which is generated by $e_j,
\omega_j^{\pm}$ (resp. $f_j, {\omega'}_j^{\pm}$), where $1\leq i
\leq 6$. Then we have

\medskip
\noindent\textbf{Proposition 1.2.} \textit{The algebra
$U_{r,s}(\mathfrak{g})$ is a Hopf algebra under the
comultiplication, the counit and the antipode below}
\begin{gather*}\Delta(\omega_i^{\pm 1})=\omega_i^{\pm 1}\otimes \omega_i^{\pm
1}, \quad
 \Delta({\omega'}_i^{\pm 1})={\omega'}_i^{\pm 1}\otimes {\omega'}_i^{\pm
 1},\\
\Delta(e_i)=e_i\otimes 1 + \omega_i \otimes e_i, \quad \Delta(f_i)=1
\otimes f_i+ f_i \otimes {\omega'}_i, \\
\varepsilon (\omega_i^{\pm 1})=\varepsilon ({\omega'}_i^{\pm 1})=1 ,
\quad \varepsilon(e_i)=\varepsilon(f_i)=0,\\
S(\omega_i^{\pm 1})=\omega_i^{\mp 1}, \quad S({\omega'}_i^{\pm
1})={\omega'}_i^{\mp 1},\\
S(e_i)=-\omega_i^{-1}e_i , \quad S(f_i)=-f_i{\omega'}_i^{-1}.
\end{gather*}

We can define the left-adjoint and the right-adjoint action in Hopf
algebra $U_{r,s}(\mathfrak{g})$ as follows
$$ad_\ell \,a(b)=\sum_{(a)}a_{(1)}b\,S(a_{(2)}), \qquad ad_ra(b)=\sum_{(a)}S(a_{(1)})b\,a_{(2)},$$
where $\Delta(a)=\sum_{(a)}a_{(1)}\otimes a_{(2)}$, $a,\, b \in
U_{r,s}(\mathfrak{g})$.

\medskip
Let $U_{r,s}(\mathfrak{n})$ (resp. $U_{r,s}(\mathfrak{n}^-)$) denote
the subalgebra of ${\mathcal{B}}$ (resp. ${\mathcal{B}'}$) generated
by $e_i$ (resp. $f_i$) for all $1 \leq i \leq 6$. Let
\begin{gather*}
U^0={\mathbb{K}}[\omega_1^{\pm},\cdots,\omega_6^{\pm},
{\omega'}_1^{\pm},\cdots,{\omega'}_6^{\pm}],\\
U_0={\mathbb{K}}[\omega_{1}^{\pm 1}, \cdots,\omega_6^{\pm}] , \qquad
U'_0={\mathbb{K}}[{\omega'}_{1}^{\pm 1}, \cdots,{\omega'}_6^{\pm}],
\end{gather*}
denote the respective Laurent polynomial subalgebras of
$U_{r,s}({\mathfrak{g}})$, ${\mathcal{B}}$ and ${\mathcal{B}}'$.
Then we have ${\mathcal{B}}=U_{r,s}({\mathfrak{n}})\rtimes U_0$, and
${\mathcal{B}}'=U'_{0}\ltimes U_{r,s}({\mathfrak{n}}^{-})$.

\medskip
Similar to the type $A$ case (see [BW1]), we have
\medskip

 \noindent\textbf{Proposition
1.3.} \textit{There exists a unique skew-dual pairing $\langle \,,
\rangle:$ ${\mathcal{B}'}({\mathfrak{g}})\times {\mathcal{B}} (
{\mathfrak{g}} ) \longrightarrow {\mathbb{Q}}(r,s)$
 of the Hopf algebras $B(\mathfrak{g})$ and $B'(\mathfrak{g})$ such that
$$\langle f_i,e_j\rangle=\delta_{i j}\frac{1}{s-r}, \eqno{(1.1)}$$
$$\langle \omega'_i,\omega_j\rangle=r^{p_{ji}} (s^{-1})^{q_{ji}},\eqno{(1.2)}$$
$$\langle {\omega'}_i^\pm,{\omega_j}^{-1}\rangle=\langle {\omega'}_i^\pm,\omega_j\rangle^{-1}=
\langle \omega'_i,\omega_j\rangle^{\mp1}, \eqno{(1.3)}$$ and all
other pairs of generators are $0$. Moreover, we have $\langle
S(a),S(b)\rangle = \langle a,b\rangle$ for $a \in {\mathcal{B}'},\,
b \in {\mathcal{B}}$.}

\medskip

 As a result of Proposition 1.3, we can display the structural
constants for type $E_6$ by a matrix $A=({\tilde a}_{ij})$, where
${\tilde a}_{ij}=\langle {\omega'}_i, \omega_j \rangle$,
$$
A=\left(\begin{array}{cccccc}
    rs^{-1} & 1 & r^{-1} & 1& 1 & 1 \\
    1 & rs^{-1} & 1 & r^{-1}& 1& 1\\
    s& 1 &rs^{-1}& r^{-1} & 1 & 1 \\
    1 & s & s &rs^{-1}  & r^{-1} & 1 \\
    1 & 1 &  1& s &rs^{-1} & r^{-1} \\
    1 & 1 &  1& 1& s & rs^{-1}
    \end{array}\right).
   $$

\medskip

\noindent\textbf{Proposition 1.4.} ([BGH1, Coro. 2.7]) \textit{For
$\zeta=\sum_{i=1}^{6}\zeta_i\alpha_i \in Q$, the defining relations
$(E2)$ in $U_{r,s}({\mathfrak{g}})$ can be rewritten as the forms
below}
\begin{gather*}\omega_{\zeta}e_i\omega^{-1}_{\zeta}=\langle
{\omega'}_i ,\omega_{\zeta}\rangle e_i, \qquad
 \omega_{\zeta}f_i\omega^{-1}_{\zeta}=\langle {\omega'}_i ,\omega_{\zeta}\rangle^{-1}
 f_i,\\
{\omega'}_{\zeta}e_i{\omega'}^{-1}_{\zeta}=\langle
{\omega'}_{\zeta}, \omega_i\rangle^{-1} e_i, \qquad
 {\omega'}_{\zeta}f_i{\omega'}^{-1}_{\zeta}=\langle {\omega'}_{\zeta}, \omega_i\rangle f_i.
\end{gather*}
{\it Then $U_{r,s}({\mathfrak{g}})=\bigoplus_{\eta \in
Q}U^{\eta}_{r,s}({\mathfrak{g}})$ is Q-graded such that
\begin{eqnarray*}
U^{\eta}_{r,s}({\mathfrak{g}})&=&\left\{\; \sum
F_{\alpha}{\omega'}_{\mu}\omega_{\nu}E_{\beta} \in U \;\right|\;
\omega_{\zeta}(F_{\alpha}{\omega'}_{\mu}\omega_{\nu}E_{\beta})\omega^{-1}_{\zeta}=\langle
{\omega'}_{\beta-\alpha}, \omega_\zeta \rangle
F_{\alpha}{\omega'}_{\mu}\omega_{\nu}E_{\beta},\\
&
&{\omega'}_{\zeta}(F_{\alpha}{\omega'}_{\mu}\omega_{\nu}E_{\beta}){\omega'}^{-1}_{\zeta}=\langle
{\omega'}_\zeta, \omega_{\beta-\alpha} \rangle^{-1}
F_{\alpha}{\omega'}_{\mu}\omega_{\nu}E_{\beta},\, with \,
\beta-\alpha=\eta\; \Bigr\},
\end{eqnarray*}
where $F_\alpha\, ($resp. $E_\alpha)$ is a certain monomial $f_{i_1}
\cdots f_{i_l}\, ($resp. $e_{i_1} \cdots e_{i_m})$ such that
$\alpha_{i_1}+ \cdots + \alpha_{i_l}=\alpha$ $($resp. $\alpha_{j_1}+
\cdots + \alpha_{j_m}=\beta$$)$.}\\

\bigskip
\centerline{\textsc{2. Lyndon words and convex PBW-type basis}}
\bigskip

Thanks to the work in [LR], [K1,2] and [R2,3], there is a
combinatorial approach to constructing an ordered basis called a
convex PBW-type basis (for definition, see [R3]) for our
$U_{r,s}(\mathfrak {n})$. In this section, we will give a
description of a convex PBW-type basis of $U_{r,s}(\mathfrak{n})$
making use of Lyndon words and $(r,s)$-bracketing.

  Let $A=\{e_1,e_2,\cdots,e_6\}$ be an ordered
alphabet set and the order is defined by $e_1<e_2< \cdots <e_6$. Let
$A^*$ be the set of all words in the alphabet set $A$ and let $u<v$
denote that word $u$ is lexicographically smaller than word $v$.

\medskip
\noindent\textbf{Definition 2.1.} \textit{A word $\ell \in A^*$ is a
Lyndon word if it is lexicographically smaller than all its proper
right factors.}

\medskip
 \noindent\textbf{Definition 2.2.} \textit{Let $\ell=uv$, we call it a Lyndon decomposition if
 $u$,
$v$ are both Lyndon words and $u$ is the shortest Lyndon word
appearing as a proper left factor of $\ell$. }

\medskip
Let ${\mathbb{K}}[A^*]$ be the associative algebra of
$\mathbb{K}$-linear combinations of words $A^*$ whose product is
juxtaposition, namely, a free ${\mathbb{K}}$-algebra.

\medskip
 \noindent\textbf{Theorem 2.3.}([LR], [R2,3]) \textit{The set of products $\ell_1 \cdots \ell_k$ is a
basis of ${\mathbb{K}}[A^*]$, where the $\ell_i$'s are Lyndon words
and $\ell_1\geq \cdots \geq \ell_k$}.

\medskip
Let $J$ be the $(r,s)$-Serre ideal of ${\mathbb{K}}[A^*]$ generated
by elements $\{(ad_\ell e_i)^{1-a_{ij}}(e_j)\mid 1\leq i\ne j
\leq6\}$. Now it is clear that
$U_{r,s}({\mathfrak{n}})={\mathbb{K}}[A^*]/J$.

In order to construct a monomial basis of $U_{r,s}(\mathfrak{n})$,
we need to give another kind of order $\preceq$ in $A^*$ with
introducing a usual length function $|\cdot|$ for a word $u\in A^*$.
We say $u\preceq w $, if $|u|<|w|$ or $|u|=|w|$ and $u\geq w$.

\medskip
\noindent\textbf{Definition 2.4.} \textit{Call a (Lyndon) word to be
good w.r.t. the $(r,s)$-Serre ideal $J$ if it cannot be written as a
sum of strictly smaller words modulo $J$ w.r.t. the ordering
$\preceq$.}

\medskip
For example, $e_1e_2$ is not ``good", since $e_1e_2=e_2e_1$ and
$e_2e_1$ is strictly smaller than $e_1e_2$ w.r.t. to the ordering
$\preceq$.

\medskip
\noindent\textbf{Theorem 2.5.} \textit{The set of products $\ell_1
\cdots \ell_k$, where $\ell_i$'s are good Lyndon words and $\ell_1
\geq\cdots \geq \ell_k$, is a basis of $U_{r,s}({\mathfrak{n}})$
$($Set $U^+:=U_{r,s}({\mathfrak{n}})$ for\, short $)$.}

\medskip
\noindent{\it Proof.} First, we claim that the set of good words is
a basis for $U_{r,s}({\mathfrak{n}})={\mathbb{K}}[A^*] / J$. Every
element in ${\mathbb{K}}[A^*] /J$ can be written as a linear
combination of the words in ${\mathbb{K}}[A^*]$ and if any of them
is not ``good", then we can change it into good ones w.r.t. to $J$.
This process can be continued until all the monomials appearing in
the linear combination are good, then we get our claim. Second, any
factor of a good word is a good word. Otherwise, if $u=u_1u_2 \cdots
u_n$ is a good word but a factor of it, say $u_i$, is not good, then
we have that $ u_i=\sum_{m\prec u_i}a_m m \quad (mod\, J)$ such that
$$u=u_1 \cdots u_{i-1}(\sum_{m\prec u_i}a_m m)u_{i+1} \cdots u_n
\quad (mod\, J).$$
That means $u$ is not a good word. It is a
contradiction. In view of Theorem 2.3, we get the result. $\qquad
\Box$

\medskip
More precisely, we have the following inductive construction. For
each pair of homogeneous elements $u \in U^+_{\zeta}, v \in
U^+_{\eta}$, we fix the notation
$p_{\zeta\eta}=\langle\omega'_{\eta},\omega_{\zeta}\rangle$, and
define a bilinear skew commutator named \textit{$(r,s)$-bracketing}
on the set of graded homogeneous noncommutative polynomials $u,\,v$
by the formula
$$\lceil u,v \rfloor=uv-p_{\zeta\eta}vu=uv-\langle\omega'_{\eta},\omega_{\zeta}\rangle vu.$$
We call $\lceil u \rfloor$ a \textit{good letter} (or say, a
\textit{quantum root vector}) in $U^+$ if $u$ is a good Lyndon word.
By induction, we define $\lceil u \rfloor$ as
$$\lceil u \rfloor=\lceil\lceil v \rfloor \lceil w \rfloor\rfloor, \quad \text{if } \ u=vw \ \text{ is a Lyndon decomposition.}$$

We list all the good Lyndon words ordered by $<$ and the figure of
them as follows\\

\hspace{12pt}

%TeXCAD Options
%\grade{\on}
%\emlines{\off}
%\epic{\off}
%\beziermacro{\on}
%\reduce{\on}
%\snapping{\off}
%\quality{8.00}
%\graddiff{0.01}
%\snapasp{1}
%\zoom{4.0000}
\unitlength 1mm % = 2.85pt
\linethickness{0.4pt}
\ifx\plotpoint\undefined\newsavebox{\plotpoint}\fi % GNUPLOT compatibility
\begin{picture}(98.5,106.25)(0,0)
\put(10.79,36.5){\line(1,0){42}} \put(10.79,36.5){\circle*{1.58}}
\put(52.79,36.5){\circle*{1.58}} \put(25.0,36.5){\circle*{1.58}}
\put(38.29,36.5){\circle*{1.58}} \put(10.04,31.5){3}
\put(24.5,31.5){4} \put(38.04,31.5){5} \put(52.5,31.5){6}
\put(10.79,26.5){\line(1,0){27.75}}
\put(10.79,26.75){\circle*{1.58}} \put(24.04,26.5){\circle*{1.58}}
\put(38.04,26.25){\circle*{1.58}} \put(10.0,21.5){4}
\put(23.75,21.5){5} \put(38.0,21.5){6}
%\emline(11.04,60.25)(19.04,53)
\multiput(11.04,60.25)(.0372093,-.03372093){215}{\line(1,0){.0372093}}
%\end
\put(11.04,60.75){\circle*{1.58}} \put(18.79,53.25){\circle*{1.58}}
%\emline(19.29,52.75)(28.04,45.25)
\multiput(19.29,52.75)(.03923767,-.03363229){223}{\line(1,0){.03923767}}
%\end
\put(28.04,45){\circle*{1.58}}
%\emline(19.54,54)(27.29,60.75)
\multiput(19.54,54)(.03855721,.03358209){201}{\line(1,0){.03855721}}
%\end
\put(27.29,60.25){\circle*{1.58}}
%\emline(28.29,60)(35.29,54)
\multiput(28.29,60)(.03932584,-.03370787){178}{\line(1,0){.03932584}}
%\end
\put(35.54,53.75){\circle*{1.58}}
%\emline(36.04,53.25)(45.04,44.75)
\multiput(36.04,53.25)(.035714286,-.033730159){252}{\line(1,0){.035714286}}
%\end
\put(44.79,44.75){\circle*{1.58}}
%\emline(27.54,61)(38.54,69.5)
\multiput(27.54,61)(.043650794,.033730159){252}{\line(1,0){.043650794}}
%\end
\put(37.79,69){\circle*{1.58}} \put(37.79,68.5){\line(6,-5){9}}
\put(47.04,60.75){\circle*{1.58}}
%\emline(47.79,60)(56.04,52.25)
\multiput(47.79,60)(.03586957,-.03369565){230}{\line(1,0){.03586957}}
%\end
\put(56.29,52.25){\circle*{1.58}}
%\emline(56.79,51.75)(65.29,44)
\multiput(56.79,51.75)(.03695652,-.03369565){230}{\line(1,0){.03695652}}
%\end
\put(65.29,44){\circle*{1.58}} \put(10.54,65.5){2} \put(18.29,58){4}
\put(25.79,65.25){5} \put(23.29,44.5){3} \put(32.29,52){3}
\put(39.54,43){4} \put(50.54,63.75){3} \put(59.54,56){4}
\put(69.29,47){5} \put(37.54,72.25){6}
%\emline(22.75,95.25)(31.5,87.5)
\multiput(22.75,95.25)(.03804348,-.03369565){230}{\line(1,0){.03804348}}
%\end
%\emline(31.5,87.5)(38.5,80.75)
\multiput(31.5,87.5)(.03482587,-.03358209){201}{\line(1,0){.03482587}}
%\end
%\emline(30.75,88.5)(38.5,95.25)
\multiput(30.75,88.5)(.03855721,.03358209){201}{\line(1,0){.03855721}}
%\end
\put(30.25,88.5){\circle*{1.58}} \put(47.25,88.75){\circle*{1.58}}
\put(38.25,81.25){\circle*{1.58}}
%\emline(38.75,95.5)(46.75,89)
\multiput(38.75,95.5)(.04145078,-.03367876){193}{\line(1,0){.04145078}}
%\end
\put(47.75,88.25){\line(1,-1){6.75}} \put(54.5,81.5){\circle*{1.58}}
%\emline(54.75,81.25)(62,74.25)
\multiput(54.75,81.25)(.03485577,-.03365385){208}{\line(1,0){.03485577}}
%\end
\put(62,74){\circle*{1.58}} \put(38.75,95.25){\circle*{1.58}}
%\emline(39.5,95.5)(49.25,103.5)
\multiput(39.5,95.5)(.04096639,.03361345){238}{\line(1,0){.04096639}}
%\end
\put(49,103.25){\circle*{1.58}}
%\emline(49.75,103)(56.75,96.5)
\multiput(49.75,103)(.03626943,-.03367876){193}{\line(1,0){.03626943}}
%\end
\put(57.25,96){\circle*{1.58}} \put(57.75,95.5){\line(1,-1){7}}
\put(64.75,88.5){\circle*{1.58}} \put(65.25,88){\line(1,-1){6}}
\put(71.5,81.5){\circle*{1.58}}
%\emline(65.25,89.25)(74,97)
\multiput(65.25,89.25)(.03804348,.03369565){230}{\line(1,0){.03804348}}
%\end
\put(74,96.75){\circle*{1.58}}
%\emline(74.25,96)(80.75,89.25)
\multiput(74.25,96)(.03367876,-.03497409){193}{\line(0,-1){.03497409}}
%\end
\put(81.25,89){\circle*{1.58}}
%\emline(82,88.5)(88.75,81.25)
\multiput(82,88.5)(.03358209,-.03606965){201}{\line(0,-1){.03606965}}
%\end
\put(88.5,81.5){\circle*{1.58}} \put(89.5,80.75){\line(1,-1){6.5}}
\put(96.25,74){\circle*{1.58}} \put(21.75,98.5){3} \put(30,91.75){4}
\put(37,98){5} \put(47.75,106.25){6} \put(44.75,87){2}
\put(52,80){4} \put(59.25,72.75){3} \put(55.25,94.25){2}
\put(62.25,87){4} \put(68.75,80.25){3} \put(72.75,99.75){5}
\put(89.75,83.5){4} \put(98.5,77.25){2} \put(82.75,91.5){3}
\put(35.25,79.75){2} \put(22.5,95.25){\circle*{1.58}}
%\emline(21.75,94.75)(13,88.25)
\multiput(21.75,94.75)(-.04533679,-.03367876){193}{\line(-1,0){.04533679}}
%\end
\put(13,88.5){\circle*{1.58}} \put(11.5,93.25){1}
\put(9.9,7.75){\circle*{1.58}} \put(9.5,3){6}
\put(10.29,17.75){\line(1,0){13.75}}
\put(10.29,17.75){\circle*{1.58}} \put(24.29,17.25){\circle*{1.58}}
\put(10.0,12.75){5} \put(24.0,12.75){6}
\end{picture}
\begin{eqnarray*}
&E_{1}&E_{13}\,\,\,\,E_{134}\,\,\,\,E_{1342}\,\,\,\,E_{1345}\,\,\,\,E_{13452}\,\,\,\,
E_{134524}\,\,\,\,E_{1345243}\,\,\,\,E_{13456}\,\,\,\,
E_{134562}\\
& \,\,\,&E_{1345624}\,\,\,\,E_{13456243}\,\,\,\,E_{13456245}\,\,\,\,
E_{134562453}\,\,\,\,E_{1345624534}\,\,\,\,E_{13456245342}\\
&E_{2}&E_{24}\,\,\,\,E_{243}\,\,\,\,E_{245}\,\,\,\,E_{2453}\,\,\,\,E_{24534}\,\,\,\,E_{2456}\,\,\,\,E_{24563}\\
& &E_{245634}\,\,\,\,E_{2456345}\\
&E_{3}&E_{34}\,\,\,\,E_{345}\,\,\,\,E_{3456}\hspace{8cm}\\
&E_{4}&E_{45}\,\,\,\,E_{456}\\
&E_{5}&E_{56}\\
&E_{6}&
\end{eqnarray*} where $E_{i_1 \cdots i_n}$ denotes
$e_{i_1}e_{i_2}\cdots
e_{i_n}$.\\

Denote ${\mathcal{E}}_{\beta_{1}}=\lceil E_1\rfloor,
{\mathcal{E}}_{\beta_{2}}=\lceil E_{13}\rfloor,
{\mathcal{E}}_{\beta_{3}}=\lceil E_{134}\rfloor,\cdots,
{\mathcal{E}}_{\beta_{36}}=\lceil E_{6}\rfloor$, where $\beta_i$
denotes a root in $\Phi^+$. Then we have the following theorem.

\medskip
\noindent\textbf{Theorem 2.6.} \textit{The set of products
${\mathcal{E}}_{\beta_{36}}^{n_{36}}\cdots
{\mathcal{E}}_{\beta_{2}}^{n_{2}}{\mathcal{E}}_{\beta_{1}}^{n_{1}}$
is a convex PBW-type basis of $U_{r,s}({\mathfrak{n}}^+)$, which is
a Lyndon basis with the ``convexity property" in the sense of
$[R3]$, where $n_1, \cdots,n_{36}$ are nonnegative integers.}

\medskip

The proof is similar to that in [K1].

\medskip
Similarly, we define a bilinear skew commutator on the set of graded
homogeneous noncommutative polynomials in
$U_{r,s}^-({\mathfrak{n}})$. For each pair of homogeneous elements
$u,v$ in the free algebra ${\mathbb{K}}\langle f_1,\cdots,f_6
\rangle$ and $u \in U^-_{\zeta}, v \in U^-_{\eta}$, we fix the
notation
${p'}_{\zeta\eta}=\langle\omega'_{\zeta},\omega_{\eta}\rangle^{-1}$,
$$\lceil u,v \rfloor=vu-{p'}_{\zeta\eta}uv=vu-\langle\omega'_{\zeta},\omega_{\eta}\rangle^{-1}uv.$$
We call $\lceil u \rfloor$ a \textsl{good letter} if $u$ is a good
Lyndon word. By induction, we define $\lceil u \rfloor$ as
$$\lceil u \rfloor=\lceil\lceil v \rfloor \lceil w \rfloor\rfloor,$$
where $u=vw$ is a Lyndon decomposition. Denote $f_{i_1}f_{i_2}\cdots
f_{i_n}$ by $F_{i_1 \cdots i_n}$ and set
${\mathcal{F}}_{\beta_{1}}=\lceil F_1\rfloor,\,
{\mathcal{F}}_{\beta_{2}}=\lceil F_{13}\rfloor,\,
{\mathcal{F}}_{\beta_{3}}=\lceil F_{134}\rfloor,\,\cdots,\,
{\mathcal{F}}_{\beta_{36}}=\lceil F_{6}\rfloor$. The set of products
${\mathcal{F}}_{\beta_{36}}^{n_{36}}{\mathcal{F}}_{\beta_{35}}^{n_{35}}\cdots
{\mathcal{F}}_{\beta_{1}}^{n_{1}}$ is a basis of
$U_{r,s}({\mathfrak{n}}^-)$, where $n_1, \cdots, n_{36}$ are
nonnegative integers.\\

\bigskip
\centerline{\textsc{3. Drinfeld double and universal $R$-matrix}}
\bigskip

In this section, we will give the Drinfeld double structure of the
algebra $U_{r,s}({\mathfrak{g}})$ after preparing some of Lemmas.
This structure, together with the result about the convex PBW-type
basis, will be used to construct the explicit form of the canonical
element and the universal $R$-matrix of  $U_{r,s}({\mathfrak{g}})$.

\medskip
 \noindent\textbf{Lemma 3.1.}
 \textit{$\Delta
({\mathcal{E}}_{\beta_{i}})={\mathcal{E}}_{\beta_{i}}\otimes 1 +
\omega_{{\mathcal{E}}_{\beta_{i}}}\otimes {\mathcal{E}}_{\beta_{i}}+
\sum(*){\mathcal{E}}_{\beta_{i}}^{(1)}\omega_{{\mathcal{E}}_{\beta_{i}}^{(2)}}
\otimes{\mathcal{E}}_{\beta_{i}}^{(2)}$, where
$deg({\mathcal{E}}_{\beta_{i}})=deg({\mathcal{E}}_{\beta_{i}}^{(1)})+deg({\mathcal{E}}_{\beta_{i}}^{(2)})$,
${\mathcal{E}}_{\beta_{i}}^{(1)}$ $(<{\mathcal{E}}_{\beta_{i}})$ is
a good letter, and ${\mathcal{E}}_{\beta_{i}}^{(2)}$ is a non
increasing product of good letters $($i.e., $(r,s)$-bracketing of
Lyndon words$)$, which are bigger than ${\mathcal{E}}_{\beta_{i}}$
w.r.t. the ordering $<$}.

\medskip
\noindent\textit{Proof}. We will prove it by induction. Assume that
the Lyndon decomposition of ${\mathcal{E}}_{\beta_{i}}$ is
${\mathcal{E}}_{\beta_{i}}=\lceil {\mathcal{E}}_{i1},
{\mathcal{E}}_{i2}\rfloor$ and  ${\mathcal{E}}_{i1},
{\mathcal{E}}_{i2}$ satisfying
\begin{gather*}
\Delta ({\mathcal{E}}_{i1})={\mathcal{E}}_{i1}\otimes 1 +
\omega_{{\mathcal{E}}_{i1}}\otimes {\mathcal{E}}_{i1}+
\sum(*){\mathcal{E}}_{i1}^{(1)}\omega_{{\mathcal{E}}_{i1}^{(2)}}\otimes{\mathcal{E}}_{i1}^{(2)},\\
\Delta ({\mathcal{E}}_{i2})={\mathcal{E}}_{i2}\otimes 1 +
\omega_{{\mathcal{E}}_{i2}}\otimes {\mathcal{E}}_{i2}+
\sum(*){\mathcal{E}}_{i2}^{(1)}\omega_{{\mathcal{E}}_{i2}^{(2)}}\otimes{\mathcal{E}}_{i2}^{(2)},
\end{gather*}
where ${\mathcal{E}}_{i1}^{(2)}$'s are non increasing products of
good letters ($>{\mathcal{E}}_{i1}$), ${\mathcal{E}}_{i2}^{(2)}$'s
are non increasing products of good letters ($>{\mathcal{E}}_{i2}$).
Then we have
\begin{equation*}
\begin{split}
\Delta ({\mathcal{E}}_{\beta_{i}})
&=\bigl({\mathcal{E}}_{i1}\otimes
1 + \omega_{{\mathcal{E}}_{i1}}\otimes {\mathcal{E}}_{i1}+
\sum(*){\mathcal{E}}_{i1}^{(1)}\omega_{{\mathcal{E}}_{i1}^{(2)}}\otimes{\mathcal{E}}_{i1}^{(2)}\bigr)\cdot\\
&\quad\; \bigl({\mathcal{E}}_{i2}\otimes 1 +
\omega_{{\mathcal{E}}_{i2}}\otimes {\mathcal{E}}_{i2}+
\sum(*){\mathcal{E}}_{i2}^{(1)}\omega_{{\mathcal{E}}_{i2}^{(2)}}\otimes{\mathcal{E}}_{i2}^{(2)}\bigr)\\
&\ -\langle\omega'_{{\mathcal{E}}_{i2}},
\omega_{{\mathcal{E}}_{i1}}\rangle \bigl({\mathcal{E}}_{i2}\otimes 1
+ \omega_{{\mathcal{E}}_{i2}}\otimes {\mathcal{E}}_{i2}+
\sum(*){\mathcal{E}}_{i2}^{(1)}\omega_{{\mathcal{E}}_{i2}^{(2)}}\otimes{\mathcal{E}}_{i2}^{(2)}\bigr)\cdot\\
&\quad\;  \bigl({\mathcal{E}}_{i1}\otimes 1 +
\omega_{{\mathcal{E}}_{i1}}\otimes {\mathcal{E}}_{i1}+
\sum(*){\mathcal{E}}_{i1}^{(1)}\omega_{{\mathcal{E}}_{i1}^{(2)}}\otimes{\mathcal{E}}_{i1}^{(2)}\bigr)\\
\end{split}
\end{equation*}
\begin{equation*}
\begin{split}
&={\mathcal{E}}_{i1}{\mathcal{E}}_{i2}\otimes
1+\omega_{{\mathcal{E}}_{i1}} \omega_{{\mathcal{E}}_{i2}}\otimes
{\mathcal{E}}_{i1}{\mathcal{E}}_{i2}+
{\mathcal{E}}_{i1}\omega_{{\mathcal{E}}_{i2}}\otimes{\mathcal{E}}_{i2}\\
& +\sum(*){\mathcal{E}}_{i1} {\mathcal{E}}_{i2}^{(1)}
\omega_{{\mathcal{E}}_{i2}^{(2)}} \otimes{\mathcal{E}}_{i2}^{(2)}
+\omega_{{\mathcal{E}}_{i1}} {\mathcal{E}}_{i2}\otimes
{\mathcal{E}}_{i1}+ \sum(*) \omega_{{\mathcal{E}}_{i1}}
{\mathcal{E}}_{i2}^{(1)} \omega_{{\mathcal{E}}_{i2}^{(2)}}
\otimes{\mathcal{E}}_{i1}{\mathcal{E}}_{i2}^{(2)}\\
 &+ \sum(*)
{\mathcal{E}}_{i1}^{(1)} \omega_{{\mathcal{E}}_{i1}^{(2)}}
{\mathcal{E}}_{i2} \otimes{\mathcal{E}}_{i1}^{(2)}+ \sum(*)
{\mathcal{E}}_{i1}^{(1)}
\omega_{{\mathcal{E}}_{i1}^{(2)}}\omega_{{\mathcal{E}}_{i2}}
\otimes{\mathcal{E}}_{i1}^{(2)} {\mathcal{E}}_{i2}\\
&+\sum(*){\mathcal{E}}_{i1}^{(1)}\omega_{{\mathcal{E}}_{i1}^{(2)}}
{\mathcal{E}}_{i2}^{(1)} \omega_{{\mathcal{E}}_{i2}^{(2)}}
\otimes{\mathcal{E}}_{i1}^{(2)}{\mathcal{E}}_{i2}^{(2)}\\ &
-\langle\omega'_{{\mathcal{E}}_{i2}},\omega_{{\mathcal{E}}_{i1}}\rangle
\bigl({\mathcal{E}}_{i2}{\mathcal{E}}_{i1}\otimes
1+\omega_{{\mathcal{E}}_{i1}} \omega_{{\mathcal{E}}_{i2}}\otimes
{\mathcal{E}}_{i2}{\mathcal{E}}_{i1}+
\omega_{{\mathcal{E}}_{i2}}{\mathcal{E}}_{i1}\otimes{\mathcal{E}}_{i2}\\
&+\sum(*){\mathcal{E}}_{i2}^{(1)} \omega_{{\mathcal{E}}_{i2}^{(2)}}
{\mathcal{E}}_{i1} \otimes{\mathcal{E}}_{i2}^{(2)} +
{\mathcal{E}}_{i2}\omega_{{\mathcal{E}}_{i1}}\otimes
{\mathcal{E}}_{i1}+ \sum(*) {\mathcal{E}}_{i2}^{(1)}
 \omega_{{\mathcal{E}}_{i2}^{(2)}}\omega_{{\mathcal{E}}_{i1}}
\otimes {\mathcal{E}}_{i2}^{(2)}{\mathcal{E}}_{i1}\\
&+ \sum(*) {\mathcal{E}}_{i2}{\mathcal{E}}_{i1}^{(1)}
\omega_{{\mathcal{E}}_{i1}^{(2)}} \otimes{\mathcal{E}}_{i1}^{(2)}+
\sum(*) \omega_{{\mathcal{E}}_{i2}}{\mathcal{E}}_{i1}^{(1)}
\omega_{{\mathcal{E}}_{i1}^{(2)}}
\otimes {\mathcal{E}}_{i2}{\mathcal{E}}_{i1}^{(2)}\\
&
+\sum(*){\mathcal{E}}_{i2}^{(1)}\omega_{{\mathcal{E}}_{i2}^{(2)}}{\mathcal{E}}_{i1}^{(1)}\omega_{{\mathcal{E}}_{i1}^{(2)}}
\otimes{\mathcal{E}}_{i2}^{(2)}{\mathcal{E}}_{i1}^{(2)}\bigr)\\
&=({\mathcal{E}}_{i1}{\mathcal{E}}_{i2}\otimes
1{-}\langle\omega'_{{\mathcal{E}}_{i2}},
\omega_{{\mathcal{E}}_{i1}}\rangle
{\mathcal{E}}_{i2}{\mathcal{E}}_{i1}\otimes 1)
+(\omega_{{\mathcal{E}}_{i1}} \omega_{{\mathcal{E}}_{i2}}\otimes
{\mathcal{E}}_{i1}{\mathcal{E}}_{i2}
{-}\langle\omega'_{{\mathcal{E}}_{i2}},
\omega_{{\mathcal{E}}_{i1}}\rangle \omega_{{\mathcal{E}}_{i1}}
\omega_{{\mathcal{E}}_{i2}}\otimes
{\mathcal{E}}_{i2}{\mathcal{E}}_{i1})\\
&+
({\mathcal{E}}_{i1}\omega_{{\mathcal{E}}_{i2}}\otimes{\mathcal{E}}_{i2}-
\langle\omega'_{{\mathcal{E}}_{i2}},
\omega_{{\mathcal{E}}_{i1}}\rangle
\omega_{{\mathcal{E}}_{i2}}{\mathcal{E}}_{i1}\otimes{\mathcal{E}}_{i2})
+(\sum(*){\mathcal{E}}_{i1} {\mathcal{E}}_{i2}^{(1)}
\omega_{{\mathcal{E}}_{i2}^{(2)}} \otimes{\mathcal{E}}_{i2}^{(2)}
-\langle\omega'_{{\mathcal{E}}_{i2}},
\omega_{{\mathcal{E}}_{i1}}\rangle \\
&\cdot\sum(*){\mathcal{E}}_{i2}^{(1)}
\omega_{{\mathcal{E}}_{i2}^{(2)}} {\mathcal{E}}_{i1}
\otimes{\mathcal{E}}_{i2}^{(2)}) +(\omega_{{\mathcal{E}}_{i1}}
{\mathcal{E}}_{i2}\otimes {\mathcal{E}}_{i1}
-\langle\omega'_{{\mathcal{E}}_{i2}},
\omega_{{\mathcal{E}}_{i1}}\rangle
{\mathcal{E}}_{i2}\omega_{{\mathcal{E}}_{i1}}\otimes
{\mathcal{E}}_{i1})\\
& + (\sum(*) \omega_{{\mathcal{E}}_{i1}} {\mathcal{E}}_{i2}^{(1)}
\omega_{{\mathcal{E}}_{i2}^{(2)}}
\otimes{\mathcal{E}}_{i1}{\mathcal{E}}_{i2}^{(2)}
-\langle\omega'_{{\mathcal{E}}_{i2}},
\omega_{{\mathcal{E}}_{i1}}\rangle \sum(*) {\mathcal{E}}_{i2}^{(1)}
 \omega_{{\mathcal{E}}_{i2}^{(2)}}\omega_{{\mathcal{E}}_{i1}}
\otimes {\mathcal{E}}_{i2}^{(2)}{\mathcal{E}}_{i1})\\
&+ (\sum(*) {\mathcal{E}}_{i1}^{(1)}
\omega_{{\mathcal{E}}_{i1}^{(2)}} {\mathcal{E}}_{i2}
\otimes{\mathcal{E}}_{i1}^{(2)}
-\langle\omega'_{{\mathcal{E}}_{i2}},
\omega_{{\mathcal{E}}_{i1}}\rangle \sum(*)
{\mathcal{E}}_{i2}{\mathcal{E}}_{i1}^{(1)}
\omega_{{\mathcal{E}}_{i1}^{(2)}} \otimes{\mathcal{E}}_{i1}^{(2)}) \\
&+ (\sum(*) {\mathcal{E}}_{i1}^{(1)}
\omega_{{\mathcal{E}}_{i1}^{(2)}}\omega_{{\mathcal{E}}_{i2}}
\otimes{\mathcal{E}}_{i1}^{(2)} {\mathcal{E}}_{i2}
-\langle\omega'_{{\mathcal{E}}_{i2}},
\omega_{{\mathcal{E}}_{i1}}\rangle \sum(*)
\omega_{{\mathcal{E}}_{i2}}{\mathcal{E}}_{i1}^{(1)}
\omega_{{\mathcal{E}}_{i1}^{(2)}}
\otimes {\mathcal{E}}_{i2}{\mathcal{E}}_{i1}^{(2)})\\
&+(\sum(*){\mathcal{E}}_{i1}^{(1)} \omega_{{\mathcal{E}}_{i1}^{(2)}}
{\mathcal{E}}_{i2}^{(1)} \omega_{{\mathcal{E}}_{i2}^{(2)}}
\otimes{\mathcal{E}}_{i1}^{(2)}{\mathcal{E}}_{i2}^{(2)}
-\langle\omega'_{{\mathcal{E}}_{i2}},
\omega_{{\mathcal{E}}_{i1}}\rangle
\sum(*){\mathcal{E}}_{i2}^{(1)}\omega_{{\mathcal{E}}_{i2}^{(2)}}
{\mathcal{E}}_{i1}^{(1)}\omega_{{\mathcal{E}}_{i1}^{(2)}}
\otimes{\mathcal{E}}_{i2}^{(2)}{\mathcal{E}}_{i1}^{(2)})\\
&=({\mathcal{E}}_{i1}{\mathcal{E}}_{i2}\otimes
1{-}\langle\omega'_{{\mathcal{E}}_{i2}},
\omega_{{\mathcal{E}}_{i1}}\rangle
{\mathcal{E}}_{i2}{\mathcal{E}}_{i1}\otimes 1)
+(\omega_{{\mathcal{E}}_{i1}} \omega_{{\mathcal{E}}_{i2}}\otimes
{\mathcal{E}}_{i1}{\mathcal{E}}_{i2}
{-}\langle\omega'_{{\mathcal{E}}_{i2}},
\omega_{{\mathcal{E}}_{i1}}\rangle \omega_{{\mathcal{E}}_{i1}}
\omega_{{\mathcal{E}}_{i2}}\otimes
{\mathcal{E}}_{i2}{\mathcal{E}}_{i1})\\
&+
({\mathcal{E}}_{i1}\omega_{{\mathcal{E}}_{i2}}\otimes{\mathcal{E}}_{i2}-
\langle\omega'_{{\mathcal{E}}_{i2}},
\omega_{{\mathcal{E}}_{i1}}\rangle
\langle\omega'_{{\mathcal{E}}_{i1}},
\omega_{{\mathcal{E}}_{i2}}\rangle
{\mathcal{E}}_{i1}\omega_{{\mathcal{E}}_{i2}}\otimes{\mathcal{E}}_{i2})\\
& +(\sum(*){\mathcal{E}}_{i1} {\mathcal{E}}_{i2}^{(1)}
\omega_{{\mathcal{E}}_{i2}^{(2)}} \otimes{\mathcal{E}}_{i2}^{(2)}
-\langle\omega'_{{\mathcal{E}}_{i2}},
\omega_{{\mathcal{E}}_{i1}}\rangle
\sum(*)\langle\omega'_{{\mathcal{E}}_{i1}},
\omega_{{\mathcal{E}}_{i2}^{(2)}}\rangle{\mathcal{E}}_{i2}^{(1)}
{\mathcal{E}}_{i1}\omega_{{\mathcal{E}}_{i2}^{(2)}}
\otimes{\mathcal{E}}_{i2}^{(2)})\\
& + (\langle \omega'_{{\mathcal{E}}_{i2}},
\omega_{{\mathcal{E}}_{i1}}\rangle
{\mathcal{E}}_{i2}\omega_{{\mathcal{E}}_{i1}} \otimes
{\mathcal{E}}_{i1}-\langle\omega'_{{\mathcal{E}}_{i2}},
\omega_{{\mathcal{E}}_{i1}}\rangle
{\mathcal{E}}_{i2}\omega_{{\mathcal{E}}_{i1}}\otimes
{\mathcal{E}}_{i1}) \\
&+ (\sum(*) \langle \omega'_{{\mathcal{E}}_{i2}^{(1)}},
\omega_{{\mathcal{E}}_{i1}}\rangle {\mathcal{E}}_{i2}^{(1)}
\omega_{{\mathcal{E}}_{i1}}\omega_{{\mathcal{E}}_{i2}^{(2)}}
\otimes{\mathcal{E}}_{i1}{\mathcal{E}}_{i2}^{(2)}
-\langle\omega'_{{\mathcal{E}}_{i2}},
\omega_{{\mathcal{E}}_{i1}}\rangle \sum(*) {\mathcal{E}}_{i2}^{(1)}
 \omega_{{\mathcal{E}}_{i2}^{(2)}}\omega_{{\mathcal{E}}_{i1}}
\otimes {\mathcal{E}}_{i2}^{(2)}{\mathcal{E}}_{i1})\\
&+ (\sum(*)\langle\omega'_{{\mathcal{E}}_{i2}},
\omega_{{\mathcal{E}}_{i1}^{(2)}}\rangle {\mathcal{E}}_{i1}^{(1)}
{\mathcal{E}}_{i2} \omega_{{\mathcal{E}}_{i1}^{(2)}}
\otimes{\mathcal{E}}_{i1}^{(2)}
-\langle\omega'_{{\mathcal{E}}_{i2}},
\omega_{{\mathcal{E}}_{i1}}\rangle \sum(*)
{\mathcal{E}}_{i2}{\mathcal{E}}_{i1}^{(1)}
\omega_{{\mathcal{E}}_{i1}^{(2)}} \otimes{\mathcal{E}}_{i1}^{(2)}) \\
&+ (\sum(*) {\mathcal{E}}_{i1}^{(1)}
\omega_{{\mathcal{E}}_{i1}^{(2)}}\omega_{{\mathcal{E}}_{i2}}
\otimes{\mathcal{E}}_{i1}^{(2)} {\mathcal{E}}_{i2}
-\langle\omega'_{{\mathcal{E}}_{i2}},
\omega_{{\mathcal{E}}_{i1}}\rangle \sum(*)
\langle\omega'_{{\mathcal{E}}_{i1}^{(1)}},
\omega_{{\mathcal{E}}_{i2}}\rangle
{\mathcal{E}}_{i1}^{(1)}\omega_{{\mathcal{E}}_{i2}}
\omega_{{\mathcal{E}}_{i1}^{(2)}}
\otimes {\mathcal{E}}_{i2}{\mathcal{E}}_{i1}^{(2)})\\
&+(\sum(*)\langle\omega'_{{\mathcal{E}}_{i2}^{(1)}},
\omega_{{\mathcal{E}}_{i1}^{(2)}} \rangle {\mathcal{E}}_{i1}^{(1)}
{\mathcal{E}}_{i2}^{(1)} \omega_{{\mathcal{E}}_{i1}^{(2)}}
\omega_{{\mathcal{E}}_{i2}^{(2)}}
\otimes{\mathcal{E}}_{i1}^{(2)}{\mathcal{E}}_{i2}^{(2)}
\\
& -\langle\omega'_{{\mathcal{E}}_{i2}},
\omega_{{\mathcal{E}}_{i1}}\rangle\sum(*)\langle\omega'_{{\mathcal{E}}_{i1}^{(1)}},
\omega_{{\mathcal{E}}_{i2}^{(2)}}\rangle {\mathcal{E}}_{i2}^{(1)}
{\mathcal{E}}_{i1}^{(1)}\omega_{{\mathcal{E}}_{i2}^{(2)}}
\omega_{{\mathcal{E}}_{i1}^{(2)}}
\otimes{\mathcal{E}}_{i2}^{(2)}{\mathcal{E}}_{i1}^{(2)})\\
&={\mathcal{E}}_{\beta_{i}}\otimes 1
+\omega_{{\mathcal{E}}_{\beta_{i}}}\otimes {\mathcal{E}}_{\beta_{i}}
+(1-r^{-1}s){\mathcal{E}}_{i1}\omega_{{\mathcal{E}}_{i2}}\otimes{\mathcal{E}}_{i2}\\
&+(\sum(*){\mathcal{E}}_{i1} {\mathcal{E}}_{i2}^{(1)}
\omega_{{\mathcal{E}}_{i2}^{(2)}} \otimes{\mathcal{E}}_{i2}^{(2)}
-\langle\omega'_{{\mathcal{E}}_{i2}},
\omega_{{\mathcal{E}}_{i1}}\rangle
\sum(*)\langle\omega'_{{\mathcal{E}}_{i1}},
\omega_{{\mathcal{E}}_{i2}^{(2)}}\rangle{\mathcal{E}}_{i2}^{(1)}
{\mathcal{E}}_{i1}\omega_{{\mathcal{E}}_{i2}^{(2)}}
\otimes{\mathcal{E}}_{i2}^{(2)}) + 0 \\
&+ (\sum(*) \langle \omega'_{{\mathcal{E}}_{i2}^{(1)}},
\omega_{{\mathcal{E}}_{i1}}\rangle {\mathcal{E}}_{i2}^{(1)}
\omega_{{\mathcal{E}}_{i1}}\omega_{{\mathcal{E}}_{i2}^{(2)}}
\otimes{\mathcal{E}}_{i1}{\mathcal{E}}_{i2}^{(2)}
-\langle\omega'_{{\mathcal{E}}_{i2}},
\omega_{{\mathcal{E}}_{i1}}\rangle \sum(*) {\mathcal{E}}_{i2}^{(1)}
 \omega_{{\mathcal{E}}_{i2}^{(2)}}\omega_{{\mathcal{E}}_{i1}}
\otimes {\mathcal{E}}_{i2}^{(2)}{\mathcal{E}}_{i1})\\
\end{split}
\end{equation*}
\begin{equation*}
\begin{split}&+ (\sum(*)\langle\omega'_{{\mathcal{E}}_{i2}},
\omega_{{\mathcal{E}}_{i1}^{(2)}}\rangle {\mathcal{E}}_{i1}^{(1)}
{\mathcal{E}}_{i2} \omega_{{\mathcal{E}}_{i1}^{(2)}}
\otimes{\mathcal{E}}_{i1}^{(2)}
-\langle\omega'_{{\mathcal{E}}_{i2}},
\omega_{{\mathcal{E}}_{i1}}\rangle \sum(*)
{\mathcal{E}}_{i2}{\mathcal{E}}_{i1}^{(1)}
\omega_{{\mathcal{E}}_{i1}^{(2)}} \otimes{\mathcal{E}}_{i1}^{(2)}) \\
&+ (\sum(*) {\mathcal{E}}_{i1}^{(1)}
\omega_{{\mathcal{E}}_{i1}^{(2)}}\omega_{{\mathcal{E}}_{i2}}
\otimes{\mathcal{E}}_{i1}^{(2)} {\mathcal{E}}_{i2}
-\langle\omega'_{{\mathcal{E}}_{i2}},
\omega_{{\mathcal{E}}_{i1}}\rangle \sum(*)
\langle\omega'_{{\mathcal{E}}_{i1}^{(1)}},
\omega_{{\mathcal{E}}_{i2}}\rangle
{\mathcal{E}}_{i1}^{(1)}\omega_{{\mathcal{E}}_{i2}}
\omega_{{\mathcal{E}}_{i1}^{(2)}}
\otimes {\mathcal{E}}_{i2}{\mathcal{E}}_{i1}^{(2)})\\
&+(\sum(*)\langle\omega'_{{\mathcal{E}}_{i2}^{(1)}},
\omega_{{\mathcal{E}}_{i1}^{(2)}} \rangle {\mathcal{E}}_{i1}^{(1)}
{\mathcal{E}}_{i2}^{(1)} \omega_{{\mathcal{E}}_{i1}^{(2)}}
\omega_{{\mathcal{E}}_{i2}^{(2)}}
\otimes{\mathcal{E}}_{i1}^{(2)}{\mathcal{E}}_{i2}^{(2)}
\\
& -\langle\omega'_{{\mathcal{E}}_{i2}},
\omega_{{\mathcal{E}}_{i1}}\rangle\sum(*)\langle\omega'_{{\mathcal{E}}_{i1}^{(1)}},
\omega_{{\mathcal{E}}_{i2}^{(2)}}\rangle {\mathcal{E}}_{i2}^{(1)}
{\mathcal{E}}_{i1}^{(1)}\omega_{{\mathcal{E}}_{i2}^{(2)}}
\omega_{{\mathcal{E}}_{i1}^{(2)}}
\otimes{\mathcal{E}}_{i2}^{(2)}{\mathcal{E}}_{i1}^{(2)})\\
&={\mathcal{E}}_{\beta_{i}}\otimes 1 +
\omega_{{\mathcal{E}}_{\beta_{i}}}\otimes {\mathcal{E}}_{\beta_{i}}+
\sum(*){\mathcal{E}}_{\beta_{i}}^{(1)}\omega_{{\mathcal{E}}_{\beta_{i}}^{(2)}}\otimes{\mathcal{E}}_{\beta_{i}}^{(2)}.
\end{split}
\end{equation*}

We can give some notes for the discussion above. Since any Lyndon
word is smaller than its proper right factors, we can deduce that
${\mathcal{E}}_{i2}, {\mathcal{E}}_{i2}^{(2)},
{\mathcal{E}}_{i1}^{(2)}, {\mathcal{E}}_{i1}^{(2)}
{\mathcal{E}}_{i2} , {\mathcal{E}}_{i2}{\mathcal{E}}_{i1}^{(2)},
{\mathcal{E}}_{i1}^{(2)}{\mathcal{E}}_{i2}^{(2)}$, and $
{\mathcal{E}}_{i2}^{(2)}{\mathcal{E}}_{i1}^{(2)}$ can be written as
non increasing products of good letters, which are bigger than
${\mathcal{E}}_{\beta_{i}}$.
%For the terms
%$$\sum(*)
%\langle \omega'_{{\mathcal{E}}_{i2}^{(1)}},
%\omega_{{\mathcal{E}}_{i1}}\rangle {\mathcal{E}}_{i2}^{(1)}
%\omega_{{\mathcal{E}}_{i1}}\omega_{{\mathcal{E}}_{i2}^{(2)}}
%\otimes{\mathcal{E}}_{i1}{\mathcal{E}}_{i2}^{(2)},$$
%we pull
%${\mathcal{E}}_{i1}$ to the right: all the terms are good except for
%$$\sum(*) \langle \omega'_{{\mathcal{E}}_{i2}^{(1)}},
%\omega_{{\mathcal{E}}_{i1}}\rangle \langle
%\omega'_{{\mathcal{E}}_{i2}^{(2)}},
%\omega_{{\mathcal{E}}_{i1}}\rangle{\mathcal{E}}_{i2}^{(1)}
%\omega_{{\mathcal{E}}_{i1}}\omega_{{\mathcal{E}}_{i2}^{(2)}}
%\otimes{\mathcal{E}}_{i1}{\mathcal{E}}_{i2}^{(2)},$$ by homogeneity
%of $\Delta$, these terms just cancel with the terms
%$$-\langle\omega'_{{\mathcal{E}}_{i2}},
%\omega_{{\mathcal{E}}_{i1}}\rangle \sum(*) {\mathcal{E}}_{i2}^{(1)}
% \omega_{{\mathcal{E}}_{i2}^{(2)}}\omega_{{\mathcal{E}}_{i1}}
%\otimes {\mathcal{E}}_{i2}^{(2)}{\mathcal{E}}_{i1},$$
So we get the lemma. \qquad $\qquad\Box$

\medskip
 Assume ${\mathcal{B}}$ is the Hopf subalgebra of $U_{r,s}({\mathfrak{g}})$ generated
 by $e_i, \omega_{i}^{\pm1} \,(1\leq i \leq 6)$, and ${\mathcal{B}'}$ is
  the Hopf subalgebra of $U_{r,s}({\mathfrak{g}})$ generated
 by $f_i, {\omega'}_{i}^{\pm1} \,(1\leq i \leq 6)$.

 Let us introduce linear forms $\eta_{\beta_{i}}$ and $\gamma_i$ in
  ${\mathcal{B}}^*$, defined by
 $$
 \eta_{\beta{i}}=\sum_{g \in G({\mathcal{B}})}({\mathcal{E}}_{\beta{i}}g)^*, \qquad
 \gamma_i(\omega_j)=\langle\omega'_i, \omega_j\rangle, \qquad \gamma_i(e_j)=0,
 $$
 where $G({\mathcal{B}})$ is the abelian group generated by $\omega_i\,(1\leq i\leq 6)$,
 and the asterisk denotes the dual basis element relative to the PBW-type basis of ${\mathcal{B}}$.
 The isomorphism $\phi : {{\mathcal{B}}'}^{coop}\rightarrow {{\mathcal{B}}}^*$ is defined by
$$\phi(\omega'_i)=\gamma_i, \qquad \phi(f_i)=\eta_i.$$
First, we will check that $\phi$ is a Hopf algebra homomorphism, and
then we will show that it is a bijection.

Now we give a series of Lemmas, with some ideas benefited from [R1].

\medskip
\noindent\textbf{Lemma 3.2.} $\gamma_i\eta_j\gamma_{i}^{-1}=\langle
\omega'_i, \omega_j \rangle \eta_j$.

\medskip
\noindent\textit{Proof.} First, we should note that $\gamma_i$'s are
invertible elements in ${\mathcal{B}}^*$ and they are commutative
with one another. It is also not difficult to see that the action of
$\gamma_i\eta_j\gamma_{i}^{-1}$ is nonzero only on basis elements of
the form $e_j\omega_{1}^{k_1}\cdots \omega_{6}^{k_6}$, and on these
elements it takes the same value
\begin{eqnarray*}
&&\gamma_i\eta_j\gamma_{i}^{-1}(e_j\omega_{1}^{k_1}\cdots
\omega_{6}^{k_6})\\
&=&\gamma_i\otimes \eta_j \otimes \gamma_{i}^{-1}((e_j\otimes
1\otimes 1 + \omega_j\otimes e_j \otimes 1 + \omega_j\otimes
\omega_j \otimes e_j)(\omega_{1}^{k_1}\cdots
\omega_{6}^{k_6})^{\otimes 3})\\
&=&\gamma_i(\omega_{j}\omega_{1}^{k_1}\cdots
\omega_{6}^{k_6})\eta_j(e_j\omega_{1}^{k_1}\cdots
\omega_{6}^{k_6})\gamma_{i}^{-1}(\omega_{1}^{k_1}\cdots
\omega_{6}^{k_6})\\
&=&\gamma_i(\omega_j)=\langle\omega'_i, \omega_j\rangle.
\end{eqnarray*}
Observing that $\eta_j(e_j\omega_{1}^{k_1}\cdots
\omega_{6}^{k_6})=1$, we have $\gamma_i\eta_j\gamma_{i}^{-1}=\langle
\omega'_i, \omega_j \rangle \eta_j$. \qquad$\qquad \Box$

\medskip
\noindent\textbf{Lemma 3.3.} \ $\Delta(\eta_{i})=\eta_{i}\otimes1
+\gamma_i\otimes \eta_{i}.\hspace{8.4cm}$

\smallskip\noindent \textit{Proof}. Since the coproduct keeps the
degree, there are only two kinds of basis elements of
${\mathcal{B}}\otimes {\mathcal{B}}$ on which $\Delta(\eta_{i})$ is
nonzero. They are $e_i\omega_{1}^{j_1}\cdots\omega_{6}^{j_6}\otimes
\omega_{1}^{k_1}\cdots\omega_{6}^{k_6}$ and
$\omega_{1}^{j_1}\cdots\omega_{6}^{j_6}\otimes e_i
\omega_{1}^{k_1}\cdots\omega_{6}^{k_6}$. Calculating the actions of
$\Delta(\eta_{i})$ on them, we get
\begin{eqnarray*}
\Delta(\eta_{i})(e_i\omega_{1}^{j_1}\cdots\omega_{6}^{j_6}\otimes
\omega_{1}^{k_1}\cdots\omega_{6}^{k_6})&=&\eta_{i}(e_i\omega_{1}^{j_1+k_1}\cdots\omega_{6}^{j_6+k_6})=1,\\
\Delta(\eta_{i})(\omega_{1}^{j_1}\cdots\omega_{6}^{j_6}\otimes e_i
\omega_{1}^{k_1}\cdots\omega_{6}^{k_6})&=&\eta_{i}(\omega_{1}^{j_1}\cdots\omega_{6}^{j_6}
e_i
\omega_{1}^{k_1}\cdots\omega_{6}^{k_6})=\langle\omega'_i,\omega_{1}^{j_1}\cdots\omega_{6}^{j_6}\rangle.
\end{eqnarray*}
Correspondingly, we have
\begin{eqnarray*}
(\eta_{i}\otimes1 +\gamma_i\otimes
\eta_{i})(e_i\omega_{1}^{j_1}\cdots\omega_{6}^{j_6}\otimes
\omega_{1}^{k_1}\cdots\omega_{6}^{k_6})&=&1,\\
(\eta_{i}\otimes1 +\gamma_i\otimes
\eta_{i})(\omega_{1}^{j_1}\cdots\omega_{6}^{j_6}\otimes e_i
\omega_{1}^{k_1}\cdots\omega_{6}^{k_6})&=&\langle\omega'_i,\omega_{1}^{j_1}\cdots\omega_{6}^{j_6}\rangle.
\end{eqnarray*}
So we get $\Delta(\eta_{i})=\eta_{i}\otimes1 +\gamma_i\otimes
\eta_{i}$. \qquad $\qquad \Box$\\

\noindent\textbf{Lemma 3.4.}
\begin{eqnarray*}
&(\text{\rm i})&\eta_{i}\eta_{j}-r^{-1}\eta_{j}\eta_{i}=
(1-r^{-1}s)\eta_{\alpha_{i}+\alpha_{j}}, \quad \text{\rm if } \  a_{ij}=-1, \ \text{\rm and }\ i<j;\\
&&\eta_{i}\eta_{j}=\eta_{j}\eta_{i},
\quad \text{ if } \ a_{ij}=0.\\
&(\text{\rm
ii})&\eta_{i}^2\eta_{j}-(r^{-1}+s^{-1})\eta_{i}\eta_{j}\eta_{i}
+(r^{-1}s^{-1})\eta_{i}\eta_{j}^2=0, \quad\text{ if } \
a_{ij}=-1, \ \text{\rm and } \  i<j;\\
&&\eta_{i}^2\eta_{j}-(r+s)\eta_{i}\eta_{j}\eta_{i}
+(rs)\eta_{i}\eta_{j}^2=0, \quad \text{\rm if } \ a_{ij}=-1, \
\text{\rm and } \ i>j.
\end{eqnarray*}

\noindent\textit{Proof}. We will give the proofs of the first
identity in (i) and the first one in (ii), and the proofs of the
others are similar. We can also ignore the $\omega_i$'s in basis
since they carry no weight on $\eta$. So in our proof we can assume
that $a_{ij}=-1 $ and $ i<j$ which implies $\alpha_{i}+\alpha_{j}
\in \Phi^+$ and $e_je_i , \lceil e_i e_j\rfloor$ are in the basis.
It is also clear from the definition of $\eta$ that
$\eta_{j}\eta_{i}(e_je_i)=1$ and zero on the other monomials,
$\eta_{\alpha_{i}+\alpha_{j}}(\lceil e_i e_j\rfloor)=1$ and zero on
the other monomials. In order to get the first identity, we need to
compute the actions of $\eta_{i}\eta_{j}$ as follows
\begin{gather*}
\eta_{i}\eta_{j}(e_je_i)= (\eta_{i}\otimes\eta_{j})( \omega_j
e_i\otimes e_j)=\langle\omega'_{i},
\omega_{j}\rangle\cdot1=r^{-1},\\
\eta_{i}\eta_{j}(\lceil e_i e_j\rfloor)= \eta_{i}\eta_{j}(e_i e_j- s
e_je_i )=(1-r^{-1}s),
\end{gather*}
where we used Lemma 3.3. Then we have
$\eta_{i}\eta_{j}-r^{-1}\eta_{j}\eta_{i}=
(1-r^{-1}s)\eta_{\alpha_{i}+\alpha_{j}}$. Left (resp. right)
multiplied by $\eta_i$ on both sides of the first identity in (i),
we get
\begin{eqnarray*}
\eta_{i}^2\eta_{j}-r^{-1}\eta_{i}\eta_{j}\eta_{i}&=&
(1-r^{-1}s)\eta_{i}\eta_{\alpha_{i}+\alpha_{j}},\\
-s^{-1}\eta_{i}\eta_{j}\eta_{i}+r^{-1}s^{-1}\eta_{j}\eta_{i}^2&=&
(1-r^{-1}s)(-s^{-1})\eta_{\alpha_{i}+\alpha_{j}}\eta_{i}.
\end{eqnarray*}
Adding the two identities together, we have
$$\eta_{i}^2\eta_{j}-(r^{-1}+s^{-1})\eta_{i}\eta_{j}\eta_{i}
+(r^{-1}s^{-1})\eta_{i}\eta_{j}^2=
(1-r^{-1}s)(\eta_{i}\eta_{\alpha_{i}+\alpha_{j}}-
s^{-1}\eta_{\alpha_{i}+\alpha_{j}}\eta_{i}).$$ Since
$2\alpha_{i}+\alpha_{j}$ is not a root, the only element in basis on
which $\eta_{\alpha_{i}+\alpha_{j}}\eta_{i}$ and
$\eta_{i}\eta_{\alpha_{i}+\alpha_{j}}$ act nontrivially is
${\mathcal{E}}_{\alpha_{i}+\alpha_{j}}e_i$. Observing that
$$\eta_{i}\eta_{\alpha_{i}+\alpha_{j}}({\mathcal{E}}_{\alpha_{i}+\alpha_{j}}e_i)=s^{-1}, \quad
\eta_{\alpha_{i}+\alpha_{j}}\eta_{i}({\mathcal{E}}_{\alpha_{i}+\alpha_{j}}e_i)=1,$$
we have
$$\eta_{i}\eta_{\alpha_{i}+\alpha_{j}}-
s^{-1}\eta_{\alpha_{i}+\alpha_{j}}\eta_{i}=0.$$ So we get the
relation
$$\eta_{i}^2\eta_{j}-(r^{-1}+s^{-1})\eta_{i}\eta_{j}\eta_{i}
+(r^{-1}s^{-1})\eta_{i}\eta_{j}^2=0. \qquad\qquad  \Box $$

\medskip
Now we want to discuss the relations between $\eta_{\beta_{i}}$ and
${\mathcal{F}}_{\beta_i}$, where $\beta_i \in \Phi^+$. We can
identify ${\mathcal{B}}^*$ with ${\mathcal{B}}'$ with the opposite
comultiplication. Let $\Delta'$ denote this opposite
comultiplication and $S'$ the antipode. With Lemmas 3.2--3.4, we
have a map ${\mathcal{B}}^*\longrightarrow {\mathcal{B}}'$
$$\eta_i\mapsto (s-r)f_i \qquad \gamma_i\mapsto \omega'_i.$$
Since this map is bijective, it is an isomorphism of Hopf algebras.

\medskip
\noindent\textbf{Definition 3.5.}([BGH1]) \textit{For any two
skew-paired Hopf algebras $\mathcal{A}$ and $\mathcal{U}$ by a
skew-dual pairing $\langle , \rangle$, one may form the Drinfel'd
double $\mathcal{D(A,U)}$, which is a Hopf algebra whose underlying
coalgebra is $\mathcal{A}\otimes \mathcal{U}$ with the tensor
product coalgebra structure and algebra structure is defined by
$$(a\otimes f)(a' \otimes f')=\sum\langle S_{\mathcal{U}}(f_{(1)}), a'_{(1)}\rangle
\langle f_{(3)}, a'_{(3)}\rangle aa'_{(2)} \otimes f_{(2)}f',$$ for
$a, a' \in \mathcal{A}$ and $f, f' \in \mathcal{U}$, and the
antipode $S$ is given by
\begin{equation*}
S(a \otimes f)=(1\otimes
S_{\mathcal{U}}(f))(S_{\mathcal{A}}(a)\otimes 1). \end{equation*}}

Similar to [BW1,3] and  [BGH1], we have

\medskip
\noindent\textbf{Theorem 3.6.} \textit{The two-parameter
quantum group $U=U_{r,s}(\mathfrak{g})$ is isomorphic to the
Drinfel'd quantum double $\mathcal{D(B,B')}$.}

\medskip
\noindent \textit{Proof}. Denote the image $e_i\otimes 1$ of $e_i$
in $\mathcal{D(B,B')}$ by $\check{e}_i$, and similarly for
$\omega_i$, $\eta_i$, and $\gamma_i$. Let
$\varphi:{\mathcal{D(B,B')}}\longrightarrow U=U_{r,s}(\mathfrak{g})$
be a map defined by:
\begin{gather*}
\varphi(\check{e}_i)=e_i ,
\qquad\varphi(\check{\eta}_i)=(s-r)f_i,\\
\varphi({\check{\omega}_i}^{\pm 1})=\omega_{i}^{\pm 1} ,
\qquad\varphi({\check{\gamma}_i}^{\pm 1})={\omega'}_{i}^{\pm 1},
\end{gather*}

From the above Lemmas, it is clear that $\varphi$ keeps the
relations in $\mathcal{B}$ and $\mathcal{B}'$. It remains to check
the mixed relations (E4). Note that
\begin{gather*}
\Delta^{(2)}(e_i)=e_i\otimes 1 \otimes 1 + \omega_i\otimes e_i\otimes 1+ \omega_i\otimes \omega_i\otimes e_i,\\
(\Delta^{(2)op})(\eta_j)=1 \otimes 1 \otimes \eta_{j}+ 1 \otimes
\eta_{j}\otimes \gamma_j+ \eta_{j}\otimes \gamma_j\otimes \gamma_j.
\end{gather*}
Using the multiplication rule in $\mathcal{D(B,B')}$, we get
$$\check{\eta}_j\check{e}_i = \delta_{i,j}(\check{\omega}_{i} + \check{e}_i\check{\eta}_j - \check{\gamma}_{j}),
\quad \textrm{\it or}\quad [\check{e}_i,
\check{\eta}_j]=\delta_{i,j}(\check{\gamma}_{i}-\check{\omega}_{i}).$$
Under $\varphi$, this corresponds to the relation
$$[e_i, (s-r)f_j]=\delta_{i,j}(\omega'_i - \omega_i),\quad \textrm{\it or}\quad [e_i, f_j]=\delta_{i,j}\frac{\omega_i - \omega'_i}{r-s},$$
which is $(E4)$.  \qquad $\qquad \Box$

\medskip
\noindent\textbf{Lemma 3.7.} \textit{Let $\beta_i$ denote a root in
 $ \Phi^+$ w.r.t. the ordering $<$, and
${\mathcal{F}}_{\beta_{i}}=\lceil\lceil{\mathcal{F}}_{\beta_{i1}}\rfloor
\lceil{\mathcal{F}}_{\beta_{i2}}\rfloor\rfloor$, then
$\eta_{\beta_i}=c_{\beta_i}{\mathcal{F}}_{\beta_i}$, where
$c_{\beta_i}$ satisfies}
$$c_{\beta_i}=-\langle {\omega'}_{\beta_{i1}}, {\omega}_{\beta_{i2}}\rangle(1-
 \langle {\omega'}_{\beta_{i2}}, {\omega}_{\beta_{i1}}\rangle\langle {\omega'}_{\beta_{i1}},
  {\omega}_{\beta_{i2}}\rangle)^{-1}c_{\beta_{i1}}c_{\beta_{i2}}.$$
\noindent\textit{Proof}. Assume that ${\eta}_{\beta_{i}}=a
{\eta}_{\beta_{i1}}{\eta}_{\beta_{i2}}+b{\eta}_{\beta_{i2}}{\eta}_{\beta_{i1}}$.
Calculating the actions of both sides of the equation on elements
${\mathcal{E}}_{\beta_{i}}$ and
${\mathcal{E}}_{\beta_{i2}}{\mathcal{E}}_{\beta_{i1}}$, we have
\begin{eqnarray*}
{\eta}_{\beta_{i}}({\mathcal{E}}_{\beta_{i}})&=&1,\\
{\eta}_{\beta_{i1}}{\eta}_{\beta_{i2}}({\mathcal{E}}_{\beta_{i}})&=&1-
 \langle {\omega'}_{\beta_{i2}}, {\omega}_{\beta_{i1}}\rangle\langle {\omega'}_{\beta_{i1}},
  {\omega}_{\beta_{i2}}\rangle,\\
{\eta}_{\beta_{i2}}{\eta}_{\beta_{i1}}({\mathcal{E}}_{\beta_{i}})&=&0,\\
{\eta}_{\beta_{i}}({\mathcal{E}}_{\beta_{i2}}{\mathcal{E}}_{\beta_{i1}})&=&0,\\
{\eta}_{\beta_{i1}}{\eta}_{\beta_{i2}}({\mathcal{E}}_{\beta_{i2}}{\mathcal{E}}_{\beta_{i1}})&=&\langle
{\omega'}_{\beta_{i1}},
  {\omega}_{\beta_{i2}}\rangle,\\
{\eta}_{\beta_{i2}}{\eta}_{\beta_{i1}}({\mathcal{E}}_{\beta_{i2}}{\mathcal{E}}_{\beta_{i1}})&=&1.
\end{eqnarray*}
Then we have
\begin{eqnarray*}
{\eta}_{\beta_{i}}&=&(1-
 \langle {\omega'}_{\beta_{i2}}, {\omega}_{\beta_{i1}}\rangle\langle {\omega'}_{\beta_{i1}},
  {\omega}_{\beta_{i2}}\rangle)^{-1}
{\eta}_{\beta_{i1}}{\eta}_{\beta_{i2}}\\
&&-\langle {\omega'}_{\beta_{i1}}, {\omega}_{\beta_{i2}}\rangle(1-
 \langle {\omega'}_{\beta_{i2}}, {\omega}_{\beta_{i1}}\rangle\langle {\omega'}_{\beta_{i1}},
  {\omega}_{\beta_{i2}}\rangle)^{-1}{\eta}_{\beta_{i2}}{\eta}_{\beta_{i1}}\\
&=&-\langle {\omega'}_{\beta_{i1}}, {\omega}_{\beta_{i2}}\rangle(1-
 \langle {\omega'}_{\beta_{i2}}, {\omega}_{\beta_{i1}}\rangle\langle {\omega'}_{\beta_{i1}},
  {\omega}_{\beta_{i2}}\rangle)^{-1}\\
  &&(c_{\beta_{i1}}c_{\beta_{i2}}
  {\mathcal{F}}_{\beta_{i2}}{\mathcal{F}}_{\beta_{i1}}-\langle {\omega'}_{\beta_{i1}}, {\omega}_{\beta_{i2}}\rangle^{-1}
  c_{\beta_{i1}}c_{\beta_{i2}}{\mathcal{F}}_{\beta_{i1}}{\mathcal{F}}_{\beta_{i2}})\\
&=&-\langle {\omega'}_{\beta_{i1}}, {\omega}_{\beta_{i2}}\rangle(1-
 \langle {\omega'}_{\beta_{i2}}, {\omega}_{\beta_{i1}}\rangle\langle {\omega'}_{\beta_{i1}},
  {\omega}_{\beta_{i2}}\rangle)^{-1}c_{\beta_{i1}}c_{\beta_{i2}}{\mathcal{F}}_{\beta_{i}}\\
&=&c_{\beta_{i}}{\mathcal{F}}_{\beta_{i}}. \qquad\qquad \Box
\end{eqnarray*}
\noindent \textbf{Lemma 3.8.}
\begin{gather*}
\langle
\eta_{\beta{i}}^{n},{\mathcal{E}}_{\beta{i}}^{n'}\rangle=\delta_{n,n'}\frac{\Psi_n(rs^{-1})}{(1-rs^{-1})^n},\tag{i}\\
\langle \eta_{\beta{36}}^{n_{36}} \cdots
\eta_{\beta{2}}^{n_{2}}\eta_{\beta{1}}^{n_{1}},\,
{\mathcal{E}}_{\beta{36}}^{n'_{36}}\cdots
{\mathcal{E}}_{\beta{2}}^{n'_{2}}{\mathcal{E}}_{\beta{1}}^{n'_{1}}\rangle
=\prod_{i=1}^{36}\delta_{n_i,n'_i}
\frac{\Psi_{n_{i}}(rs^{-1})}{(1-rs^{-1})^{n_{i}}},\tag{ii}
\end{gather*}
where $\Psi_n(a)=(1-a)(1-a^2) \cdots (1-a^n)$.

\smallskip
\noindent \textit{Proof}. (i) By Lemma 3.1, we have
$$\Delta({\mathcal{E}}_{\beta_{i}})={\mathcal{E}}_{\beta_{i}}\otimes 1 +
 \omega_{{\mathcal{E}}_{\beta_{i}}}\otimes {\mathcal{E}}_{\beta_{i}} +
 \sum(*){\mathcal{E}}_{\beta_{i}}^{(1)}\omega_{{\mathcal{E}}_{\beta_{i}}^{(2)}}
 \otimes {\mathcal{E}}_{\beta_{i}}^{(2)},$$
where
$deg({\mathcal{E}}_{\beta_{i}})=deg({\mathcal{E}}_{\beta_{i}}^{(1)})+deg({\mathcal{E}}_{\beta_{i}}^{(2)})$,
and ${\mathcal{E}}_{\beta_{i}}^{(2)}$ are the products of good
letters, which are bigger than ${\mathcal{E}}_{\beta_{i}}$. From
$\Delta({\mathcal{E}}_{\beta_{i}}^{n'})=\sum {\mathcal{E}} \omega_{
\overline{{\mathcal{E}}}}\otimes \overline{{\mathcal{E}}}$, we know
that all the $\overline{{\mathcal{E}}}$'s are bigger than
${\mathcal{E}}_{\beta_{i}}$ except for ${\mathcal{E}}_{\beta_{i}}$
itself. So the terms paired with $\eta_{\beta_{i}}^{n-1}\otimes
\eta_{\beta_{i}}$ being nonzero are $\sum{\mathcal{E}}_{\beta_{i}}
\cdots{\mathcal{E}}_{\beta_{i}}\omega_{{\mathcal{E}}_{\beta_{i}}}
{\mathcal{E}}_{\beta_{i}}\cdots{\mathcal{E}}_{\beta_{i}}\otimes
{\mathcal{E}}_{\beta_{i}}$, which gives the following
\begin{equation*}
\begin{split}
\langle \eta_{\beta_{i}}^{n},{\mathcal{E}}_{\beta_{i}}^{n'}\rangle
&=\langle \eta_{\beta_{i}}^{n-1}\otimes \eta_{\beta_{i}},
\Delta({\mathcal{E}}_{\beta_{i}}^{n'})\rangle\\
&=\langle \eta_{\beta_{i}}^{n-1}\otimes \eta_{\beta_{i}},
\sum{\mathcal{E}}_{\beta_{i}}
\cdots{\mathcal{E}}_{\beta_{i}}\omega_{{\mathcal{E}}_{\beta_{i}}}
{\mathcal{E}}_{\beta_{i}}\cdots{\mathcal{E}}_{\beta_{i}}\otimes
{\mathcal{E}}_{\beta_{i}}\rangle\\
&=(1{+}\langle \omega'_{{\mathcal{E}}_{\beta_{i}}},
\omega_{{\mathcal{E}}_{\beta_{i}}}\rangle {+} \langle
\omega'_{{\mathcal{E}}_{\beta_{i}}},
\omega_{{\mathcal{E}}_{\beta_{i}}}\rangle^2{+} \cdots {+} \langle
\omega'_{{\mathcal{E}}_{\beta_{i}}},
\omega_{{\mathcal{E}}_{\beta_{i}}}\rangle^{(n'-1)} )\langle
\eta_{\beta_{i}}^{n-1},{\mathcal{E}}_{\beta_{i}}^{n'-1}\rangle\\
&=\frac{1-\langle \omega'_{{\mathcal{E}}_{\beta_{i}}},
\omega_{{\mathcal{E}}_{\beta_{i}}}\rangle^n}{1-\langle
\omega'_{{\mathcal{E}}_{\beta_{i}}},
\omega_{{\mathcal{E}}_{\beta_{i}}}\rangle}\langle
\eta_{\beta_{i}}^{n-1},{\mathcal{E}}_{\beta_{i}}^{n'-1}\rangle\\
&=\delta_{n,n'}\frac{\Psi_n(rs^{-1})}{(1-rs^{-1})^n}.
\end{split}
\end{equation*}

(ii) Similarly, we can use Lemma 3.1 to prove (ii). Since
${\mathcal{E}}_{\beta_{1}}<{\mathcal{E}}_{\beta_{2}}< \cdots
<{\mathcal{E}}_{\beta_{36}}$, the terms paired with
$\eta_{\beta_{36}}^{n_{36}}\eta_{\beta_{35}}^{n_{35}} \cdots
\eta_{\beta_{2}}^{n_{2}}\otimes \eta_{\beta_{1}}^{n_1}$ being
nonzero are of the form $?\otimes {\mathcal{E}}_{\beta_{1}}^{n_1}$.
So we have
\begin{equation*}
\begin{split}
\langle \eta_{\beta_{36}}^{n_{36}}& \cdots
\eta_{\beta_{2}}^{n_{2}}\eta_{\beta_{1}}^{n_{1}},\,
{\mathcal{E}}_{\beta_{36}}^{n'_{36}}\cdots
{\mathcal{E}}_{\beta_{2}}^{n'_{2}}{\mathcal{E}}_{\beta_{1}}^{n'_{1}}\rangle\\
&=\langle \eta_{\beta_{36}}^{n_{36}} \cdots
\eta_{\beta_{2}}^{n_{2}}\otimes \eta_{\beta_{1}}^{n_1},\,
\Delta({\mathcal{E}}_{\beta_{36}}^{n'_{36}}\cdots{\mathcal{E}}_{\beta_{2}}^{n'_{2}})
\Delta({\mathcal{E}}_{\beta_{1}}^{n'_{1}})\rangle \\
&=\frac{\Phi_{n_{1}}(rs^{-1})}{(1-rs^{-1})^{n_{1}}}\langle
\eta_{\beta_{36}}^{n_{36}}\cdots \eta_{\beta_{2}}^{n_{2}},\,
{\mathcal{E}}_{\beta_{36}}^{n'_{36}}\cdots
{\mathcal{E}}_{\beta_{2}}^{n'_{2}}\rangle\\
&=\prod_{i=1}^{36}\delta_{n_i,n'_i}
\frac{\Psi_{n_{i}}(rs^{-1})}{(1-rs^{-1})^{n_{i}}}.
\end{split}
\end{equation*}

We complete the proofs. \qquad $\qquad\Box$\\

\noindent\textbf{Theorem 3.9.}\textit{ The canonical element $\Theta
\in U_{r,s}({\mathfrak{n}}^-)
 \otimes U_{r,s}({\mathfrak{n}}^+)$ is given by
\begin{eqnarray*}
\Theta&=&\sum\frac{(1-rs^{-1})^{n_{1}}(1-rs^{-1})^{n_{2}}\cdots(1-rs^{-1})^{n_{36}}}
{\Psi_{n_{1}}(rs^{-1})\Psi_{n_{2}}(rs^{-1})\cdots
\Psi_{n_{36}}(rs^{-1})} \eta_{\beta_{36}}^{n_{36}}
\cdots\eta_{\beta_{2}}^{n_2}\eta_{\beta_{1}}^{n_{1}}\otimes{\mathcal{E}}_{\beta_{36}}^{n_{36}}\cdots
{\mathcal{E}}_{\beta_{2}}^{n_{2}}{\mathcal{E}}_{\beta_{1}}^{n_{1}}\\
&=&\sum\frac{(1-rs^{-1})^{n_{1}+\cdots + n_{36}} c_{\beta_1}^{n_{1}}
\cdots c_{\beta_{36}}^{n_{36}}}
{\Psi_{n_{1}}(rs^{-1})\Psi_{n_{2}}(rs^{-1})\cdots
\Psi_{n_{36}}(rs^{-1})} {\mathcal{F}}_{\beta_{36}}^{n_{36}}\cdots
{\mathcal{F}}_{\beta_{2}}^{n_{2}}{\mathcal{F}}_{\beta_{1}}^{n_{1}}\otimes{\mathcal{E}}_{\beta_{36}}^{n_{36}}\cdots
{\mathcal{E}}_{\beta_{2}}^{n_{2}}{\mathcal{E}}_{\beta_{1}}^{n_{1}}.
\end{eqnarray*}}

We want to describe two linear transformations $P,\widetilde{f}$,
which build up the universal $R$-matrix $\mathcal{R}$.\\

(i) $P:M' \otimes M
 \longrightarrow M \otimes M'$ is the flip operator given by
 $P(m'\otimes m)=(m\otimes m')$.

(ii) $\widetilde{f}:M \otimes M'
 \longrightarrow M \otimes M'$ is a linear transformation based on the $f$ defined below.
 \\

 We define $f:\Lambda \times \Lambda
 \longrightarrow \mathbb{K}$ as
 $$f(\lambda , \mu)=\langle \omega_\mu' , \omega_ \lambda\rangle^{-1},$$
 which satisfies
 \begin{gather*}
 f(\lambda + \mu ,\nu)=f(\lambda ,\nu)f(\mu ,\nu),\qquad
 f(\lambda ,\mu +\nu)=f(\lambda ,\mu)f(\lambda ,\nu),\\
 f(\alpha_i , \mu)=\langle \omega_\mu' , \omega_ i\rangle^{-1},\qquad
 f(\lambda , \alpha_i)=\langle\omega_ i' ,  \omega_\lambda
 \rangle^{-1}.
 \end{gather*}
Now we define linear transformations
$\widetilde{f}=\widetilde{f}_{M,M'}:M \otimes M'
 \longrightarrow M \otimes M'$ by
  $$\widetilde{f}(m \otimes m')=f(\lambda,\mu)(m \otimes
 m')$$
for $m \in M_\lambda$ and $m' \in M'_\mu $.

\medskip
\noindent\textbf{Proposition 3.10.} \textit{Let $M$ and $M'$ be any
$U_{r,s}(\mathfrak{g})$-modules in category $\mathcal{O}$ $($see
$[BW2]$, $[BGH2])$, then the map
$${\mathcal{R}}_{M',M}=\Theta \circ \widetilde{f} \circ P :M'\otimes M \longrightarrow M\otimes M' $$
is an isomorphism of $U_{r,s}(\mathfrak{g})$-modules.}

\medskip
The proof is similar to that of Theorem 3.4 in [BGH2]. On the other
hand, it is not difficult to check that each map
${\mathcal{R}}_{M,M}$ satisfies the quantum Yang-Baxter equation and
the braid relation with a
twist.\\

\bigskip
 \centerline {\textsc{4. Weight modules of finite-dimension}}

\bigskip
 Let $\Lambda$ be the weight lattice
of $\mathfrak{g}$. Associated to any $\lambda \in \Lambda$ is an
algebra homomorphism $\hat{\lambda}$ from the subalgebra $U^0$
generated by the elements $\omega_i,\omega_i'\ (1\leq i \leq 6)$ to
$\mathbb{K}$ satisfying
$$\hat{\lambda}(\omega_i)=\langle \omega'_\lambda,\omega_i\rangle,
\qquad \hat{\lambda}(\omega_i')=\langle
\omega'_i,\omega_\lambda\rangle^{-1}.$$ Let $M$ be a module for
$U_{r,s}({\mathfrak{g}})$ of dimension $d<\infty$. As $\mathbb{K}$
is algebraically closed, by linear algebra, we have
$$M=\bigoplus_\chi M_\chi,$$
where each $\chi : U^0 \rightarrow \mathbb{K}$ is an algebra
homomorphism, and $M_\chi$ is the corresponding weight space. We say
that $U^0$ acts semisimply on $M$ if $M$ can be decomposed
into genuine eigenspaces relative to $U^0$.\\

We can deduce from the relations (E2) \& (E3) that
$$e_jM_\chi \subseteq M_{\chi \cdot \widehat{\alpha_j}},\quad f_jM_\chi \subseteq
M_{\chi \cdot \widehat{-\alpha_j}}. \eqno(4.1)$$

\noindent\textbf{Lemma 4.1.} \textit{Assume that $rs^{-1}$ is not a
root of unity, and suppose $\hat{\zeta}=\hat{\eta}$ for $\zeta,\,
\eta \in Q$, then $\zeta=\eta$.}

\medskip \noindent\textit{Proof.}  Assume $\zeta=\sum_{i=1}^6
\zeta_i\alpha_i\in\Lambda$. By definition, we have
\begin{gather*}
\hat{\zeta}(\omega_j)=\langle\omega'_{\zeta},\omega_j\rangle=
r^{\sum_{i=1}^6\zeta_i p_{ji}}s^{-\sum_{i=1}^6\zeta_i q_{ji}},\\
\hat{\zeta}(\omega'_j)=\langle\omega'_j,\omega_{\zeta}\rangle^{-1}
=r^{-\sum_{i=1}^6\zeta_ip_{ij}}s^{\sum_{i=1}^6\zeta_i q_{ij}}.
\end{gather*}
If $\eta=\sum_{i=1}^6 \eta_i\alpha_i$, then the condition
$\hat{\zeta}=\hat{\eta}$ gives the equations
\begin{gather*}
r^{\sum_{i=1}^6\zeta_ip_{ji}}s^{-\sum_{i=1}^6\zeta_iq_{ji}}=r^{\sum_{i=1}^6\eta_ip_{ji}}s^{-\sum_{i=1}^6\eta_iq_{ji}},\\
r^{-\sum_{i=1}^6\zeta_ip_{ij}}s^{\sum_{i=1}^6\zeta_iq_{ij}}=r^{-\sum_{i=1}^6\eta_ip_{ij}}s^{\sum_{i=1}^6\eta_iq_{ij}}.
\end{gather*}
It is not difficult to get the equations as follows
$$r^{\sum_{i=1}^6(\zeta_i-\eta_i)p_{ji}}s^{-\sum_{i=1}^6(\zeta_i-\eta_i)q_{ji}}=1,\eqno(4.2)$$
$$r^{\sum_{i=1}^6(\zeta_i-\eta_i)p_{ij}}s^{-\sum_{i=1}^6(\zeta_i-\eta_i)q_{ij}}=1.\eqno(4.3)$$
Multiplying (4.2) with (4.3), we get
$$r^{\sum_{i=1}^6(\zeta_i-\eta_i)(p_{ji}+p_{ij})}s^{-\sum_{i=1}^6(\zeta_i-\eta_i)(q_{ji}+q_{ij})}=1. \eqno(4.4)$$
By the definitions of $p_{ij}$ and $q_{ij}$, and Lemma 1.1, we have
$$(rs^{-1})^{\sum_{i=1}^6(\zeta_i-\eta_i)(\alpha_i , \alpha_j)}=1.$$
Due to $rs^{-1}$ being not a root of unity, we get
$$\sum_{i=1}^6(\zeta_i-\eta_i)(\alpha_i , \alpha_j)=0. \eqno(4.5)$$
Since $j$ in (4.5) is arbitrary, we get a system of homogeneous
linear equations in variables $\zeta_i-\eta_i$, whose
coefficient-matrix is exactly the Cartan matrix $A$ which is
invertible. Thus, we see that the system of homogeneous linear
equations has only zero solution, that is, for any $i$, we have
$$\zeta_i-\eta_i=0.$$
So we get the result. \qquad $\qquad \Box$

\medskip
\noindent\textbf{Remark 4.2.} Owing to Lemma 4.1, we can simplify
the notation by writing $M_{\lambda}$ (for $\lambda \in \Lambda$) as
usual for the weight space instead of $M_{\hat{\lambda}}$. So it
makes sense to let (4.1) take the classical forms as
$e_jM_{\lambda}\subseteq M_{\lambda+\alpha}$ and
$f_jM_{\lambda}\subseteq M_{\lambda-\alpha}$.

\medskip
\noindent\textbf{Proposition 4.3.} \textit{If $M$ is a
finite-dimensional $U_{r,s}(\mathfrak{g})$-module and $rs^{-1}$ is
not a root of unity, then the elements $e_i,f_i$ act nilpotently on
$M$, where $1\leq i \leq 6$.}

\medskip
\noindent\textit{Proof.} Since $M$ is a direct sum of its weight
spaces, we only need to consider the actions of $e_i,\,f_i$ on each
$M_\lambda$. We know that $e_j^k.M_\lambda \subseteq M_{\lambda -
k\alpha_j}$. Since $k\alpha_j$'s are distinct and $M$ is a
finite-dimensional $U_{r,s}(\mathfrak{g})$-module, $e_i$'s act
nilpotently on $M_\lambda$. So do the actions of $f_i$'s on $M$.
\qquad $\qquad \Box$

\medskip
It is not difficult to see that any simple
$U_{r,s}(\frak{g})$-module is a highest weight module by Proposition
4.3 and (4.1). Having Lemma 4.1 for the type $E$-series, one has a
similar weight representation theory as in [BW1] for type $A$, and
[BGH2] for types $B,\,C,\,D$.\\

\bigskip
\centerline {\textsc{5. Isomorphisms among quantum groups }}

\bigskip
In what follows, we will discuss the isomorphic relationship between
the two-parameter quantum group and the one-parameter quantum double
for type $E_6$. In fact, the following result with an analogous
argument still holds for those of types $A$ (with rank $\ge3$), $D$,
and $E_7$, $E_8$.

\medskip
\noindent\textbf{Proposition 5.1.} \textit{Assume that there is an
isomorphism of Hopf algebras $\varphi :
U_{r,s}({\mathfrak{g}})\longrightarrow U_{q,q^{-1}}(\mathfrak{g})$
for some $q$, then $r = q$ and $s = q^{-1}$.}

\medskip
\noindent\textit{Proof.} Let $\pi$ be the canonical surjection from
$U_{q,q^{-1}}(\mathfrak{g})$ onto the standard one-parameter quantum
group $U_{q}(\mathfrak{g})$ of [Ja] given by $\pi (e_i)=E_i$, $\pi
(f_i)=F_i$, $\pi ({\omega_i}^{\pm 1})=K_{i}^{\pm 1}$, $\pi
({\omega'_{i}}^{\pm 1})=K_{i}^{\mp 1}$. For $1\leq i \leq 6$, we
have
$$\Delta(\pi\varphi(e_i))=(\pi\varphi\otimes \pi\varphi)\circ(\Delta (e_i)).\eqno(5.1)$$
Note that $\pi\varphi(e_i)$ is a skew-primitive element and
$\pi\varphi(\omega_i)$ is a group-like element in
$U_{q}(\mathfrak{g})$. The elements in the group $G$ generated by
$K_i$ and the skew-primitive elements span the subspace
$$\sum_{j=1}^{6}({\mathbb{K}} E_j + {\mathbb{K}} F_j)+{\mathbb{K}} G.$$
So we can assume that
$$\pi\varphi(e_i)=\sum_{j=1}^6 a_{ij}E_j+b_{ij}F_j+\sum_{g \in G}c_{ig}g,$$
where $1\leq i \leq 6$, and $a_{ij},b_{ij},c_{ig} \in {\mathbb{K}}$.
Then we have
$$\Delta(\pi\varphi(e_i))=\sum_{j=1}^6
a_{ij}(E_j\otimes 1 + K_j\otimes E_j)+ b_{ij}(1\otimes F_j +
F_j\otimes K_j^{-1})+\sum_{g \in G}c_{ig}g\otimes g. \eqno(5.2)$$ On
the other hand, we have
$$
(\pi\varphi\otimes\pi\varphi) (\Delta (e_i))=\sum_{j=1}^6
(a_{ij}E_j\otimes 1+b_{ij}F_j\otimes1+ \pi\varphi(\omega_i)\otimes
a_{ij}E_j+\pi\varphi(\omega_i)\otimes b_{ij}F_j )\eqno(5.3)$$
$$
+\sum_g (c_{ig}g\otimes 1 +\pi\varphi(\omega_i)\otimes c_{ig}g).$$
Observing the formula (5.1) and comparing the coefficients of the
terms $?\otimes 1$ in both equations above, we have
$$\sum_{j=1}^6 (a_{ij}E_j+c_{i1}1)=\sum_{j=1}^6 (a_{ij}E_j+ b_{ij}F_j )+
\sum_{g \in G}c_{ig}g+c_{i1}\pi\varphi(\omega_i).$$ Then all
$b_{ij}=0$ and $c_{ig}=0$ for all $g$ except for $g\in \{ 1 ,
\pi\varphi(\omega_i)\}$, and in which case we have
$c_{i,\pi\varphi(\omega_i)}=-c_{i1}$. So we have
$$\pi\varphi(e_i)=\sum_{j=1}^6 a_{ij}E_j+c_{i1}(1-\pi\varphi(\omega_i)).$$
Thus we can simplify the right-hand sides of the equations (5.2),
(5.3) and get
\begin{equation*}
\begin{split}
\sum_{j=1}^6 &a_{ij}(E_j\otimes 1 + K_j\otimes E_j)+c_{i1}(1\otimes
1-\pi\varphi(\omega_i)\otimes\pi\varphi(\omega_i))\\
&=\sum_{j=1}^6 a_{ij}(E_j\otimes 1 + \pi\varphi(\omega_i)\otimes
E_j)+c_{i1}(1\otimes 1-\pi\varphi(\omega_i)\otimes1)+
c_{i1}(\pi\varphi(\omega_i)\otimes
1\\
&\qquad -\pi\varphi(\omega_i)\otimes\pi\varphi(\omega_i)).
\end{split}
\end{equation*}
This implies
$$a_{ij}(K_j - \pi\varphi(\omega_i))=0, \quad\textrm{for all }\, 1\leq i \leq6.$$
So all $a_{ij}$ equal zero except for one index $j$. That means, the
index $j$ is related to the index $i$ via $\varphi$. We thus let
$j_i$ indicate such a $j$, such that
$\pi\varphi(\omega_i)=K_{j_{i}}$.

As $\omega_{i}e_k\omega_{i}^{-1}=r^{p_{ik}} (s^{-1})^{q_{ik}}e_k$
and by the results above, we get that
\begin{equation*}
\begin{split}
\pi\varphi(\omega_i e_k)&=\pi\varphi (r^{p_{ik}}s^{-q_{ik}}e_k\omega_i),\\
K_{j_{i}}(a_{k j_{k}}E_{j_{k}}+c_{k1}(1-K_{j_{k}}))&=
r^{p_{ik}}s^{-q_{ik}}(a_{k j_{k}}E_{j_{k}}+c_{k1}(1-K_{j_{k}}))K_{j_{i}},\\
q^{\langle\alpha_{j_{i}},\alpha_{j_{k}}\rangle}a_{k
j_{k}}E_{j_{k}}K_{j_{i}}+c_{k 1}(1-K_{j_{k}})K_{j_{k}}&=
r^{p_{ik}}s^{-q_{ik}}(a_{k
j_{k}}E_{j_{k}}+c_{k1}(1-K_{j_{k}}))K_{j_{i}}.
\end{split}
\end{equation*}
The last identity implies that $c_{k1}=0$ and
$q^{\langle\alpha_{j_{i}},
\alpha_{j_{k}}\rangle}=r^{p_{ik}}s^{-q_{ik}}$.

Since $(\alpha_{j_{i}},\alpha_{j_{k}})=(\alpha_{i},\alpha_{k})$, it
is not difficult to get that $r=q $ and $s=q^{-1}$ by analyzing the
three cases
\begin{eqnarray*}
\left\{ \begin{aligned}
         q^{\langle\alpha_{j_{i}}, \alpha_{j_{k}}\rangle}&=q^2, \ r^{p_{ik}}s^{-q_{ik}}=rs^{-1},\,\,\,\,\,\textrm{if} \quad i=k,\\
     q^{\langle\alpha_{j_{i}}, \alpha_{j_{k}}\rangle}&=q^{-1}, \
     r^{p_{ik}}s^{-q_{ik}}=s,\,\,\,\,\,\quad \textrm{if} \quad  a_{ik}=-1, \ i<k,\\
 q^{\langle\alpha_{j_{i}}, \alpha_{j_{k}}\rangle}&=q^{-1}, \
     r^{p_{ik}}s^{-q_{ik}}=r^{-1},\,\,\,\,\,\textrm{if} \quad
     a_{ik}=-1, \ i>k.
        \end{aligned} \right.
 \end{eqnarray*}

So we complete the proof. \qquad $\qquad\Box$\\

\bigskip
\centerline {\textsc{6. Appendix }}

\bigskip
As an interpretation of Lemma 3.1, we give an example in the
following, where the good Lyndon words arising from the type $E_6$
case.

\smallskip
By virtue of the coproduct formula for the type $A$ case, we have
$$\Delta({\mathcal{E}}_{245})={\mathcal{E}}_{245}\otimes 1 + (1-r^{-1}s)e_{2}\omega_{45}\otimes{\mathcal{E}}_{45}
+ (1-r^{-1}s){\mathcal{E}}_{24}\omega_{5}\otimes
e_{5}+\omega_{245}\otimes {\mathcal{E}}_{245}.$$

Since ${\mathcal{E}}_{2453}=\lceil{\mathcal{E}}_{245},
e_{3}\rfloor$, we have
\begin{equation*}
\begin{split}
&\Delta({\mathcal{E}}_{2453})= \Delta (\lceil{\mathcal{E}}_{245},
e_{3}\rfloor)\\
&=({\mathcal{E}}_{245}\otimes 1 +
(1{-}r^{-1}s)e_{2}\omega_{45}\otimes{\mathcal{E}}_{45} +
(1{-}r^{-1}s){\mathcal{E}}_{24}\omega_{5}\otimes
e_{5}+\omega_{245}\otimes {\mathcal{E}}_{245})\cdot\\
&\qquad (e_{3}\otimes 1 + \omega_{3}\otimes
e_{3})-r^{-1}(e_{3}\otimes 1 + \omega_{3}\otimes
e_{3})({\mathcal{E}}_{245}\otimes 1 +
(1{-}r^{-1}s)e_{2}\omega_{45}\otimes{\mathcal{E}}_{45}\\
&\quad  + (1{-}r^{-1}s){\mathcal{E}}_{24}\omega_{5}\otimes
e_{5}+\omega_{245}\otimes {\mathcal{E}}_{245})\\
&={\mathcal{E}}_{245}e_{3}\otimes 1 +
{\mathcal{E}}_{245}\omega_{3}\otimes
e_{3}+(1{-}r^{-1}s)e_{2}\omega_{45}e_{3}\otimes{\mathcal{E}}_{45}
+(1{-}r^{-1}s)e_{2}\omega_{345}\otimes{\mathcal{E}}_{45}e_{3}\\
&\quad +(1{-}r^{-1}s){\mathcal{E}}_{24}\omega_{5}e_{3}\otimes e_{5}
+(1{-}r^{-1}s){\mathcal{E}}_{24}\omega_{35}\otimes e_{5}e_{3} +
\omega_{245}e_{3}\otimes {\mathcal{E}}_{245} \\
&\quad +\omega_{2453}\otimes
{\mathcal{E}}_{245}e_{3}-r^{-1}(e_{3}{\mathcal{E}}_{245}\otimes 1+
\omega_{3}{\mathcal{E}}_{245}\otimes e_{3}+
(1-r^{-1}s)e_{3}e_{2}\omega_{45}\otimes{\mathcal{E}}_{45}\\
&\quad + (1-r^{-1}s)\omega_{3}e_{2}\omega_{45}\otimes
e_{3}{\mathcal{E}}_{45}+(1-r^{-1}s)e_{3}{\mathcal{E}}_{24}\omega_{5}\otimes
e_{5}\\
&\quad  + (1-r^{-1}s)\omega_{3}{\mathcal{E}}_{24}\omega_{5}\otimes
e_{3}e_{5} +e_{3}\omega_{245}\otimes
{\mathcal{E}}_{245}+\omega_{2453}\otimes e_{3}{\mathcal{E}}_{245})\\
&={\mathcal{E}}_{245}e_{3}\otimes 1 +
{\mathcal{E}}_{245}\omega_{3}\otimes
e_{3}+(1{-}r^{-1}s)r^{-1}e_{2}e_{3}\omega_{45}\otimes{\mathcal{E}}_{45}
+(1{-}r^{-1}s)e_{2}\omega_{345}\otimes{\mathcal{E}}_{45}e_{3}\\
&\quad +(1{-}r^{-1}s){\mathcal{E}}_{24}e_{3}\omega_{5}\otimes e_{5}
+(1{-}r^{-1}s){\mathcal{E}}_{24}\omega_{35}\otimes e_{5}e_{3}
+r^{-1}
e_{3}\omega_{245}\otimes {\mathcal{E}}_{245} \\
&\quad +\omega_{2453}\otimes
{\mathcal{E}}_{245}e_{3}-r^{-1}(e_{3}{\mathcal{E}}_{245}\otimes 1+
s{\mathcal{E}}_{245}\omega_{3}\otimes e_{3}+
(1-r^{-1}s)e_{3}e_{2}\omega_{45}\otimes{\mathcal{E}}_{45}\\
&\quad + (1-r^{-1}s)e_{2}\omega_{3}\omega_{45}\otimes
e_{3}{\mathcal{E}}_{45}+(1-r^{-1}s)e_{3}{\mathcal{E}}_{24}\omega_{5}\otimes
e_{5}\\
&\quad  + (1-r^{-1}s)s{\mathcal{E}}_{24}\omega_{35}\otimes
e_{3}e_{5} +e_{3}\omega_{245}\otimes
{\mathcal{E}}_{245}+\omega_{2453}\otimes
e_{3}{\mathcal{E}}_{245})\\
&=({\mathcal{E}}_{245}e_{3}\otimes
1-r^{-1}e_{3}{\mathcal{E}}_{245}\otimes 1)+
({\mathcal{E}}_{245}\omega_{3}\otimes
e_{3}-r^{-1}s{\mathcal{E}}_{245}\omega_{3}\otimes e_{3})\\
&\quad
+(1-r^{-1}s)(r^{-1}e_{2}e_{3}\omega_{45}\otimes{\mathcal{E}}_{45}-
r^{-1}e_{3}e_{2}\omega_{45}\otimes{\mathcal{E}}_{45})\\
&\quad
+(1-r^{-1}s)(e_{2}\omega_{345}\otimes{\mathcal{E}}_{45}e_{3}-r^{-1}e_{2}\omega_{3}\omega_{45}\otimes
e_{3}{\mathcal{E}}_{45})\\
&\quad +(1-r^{-1}s)({\mathcal{E}}_{24}e_{3}\omega_{5}\otimes
e_{5}-r^{-1}e_{3}{\mathcal{E}}_{24}\omega_{5}\otimes e_{5})\\
&\quad  +(1-r^{-1}s)({\mathcal{E}}_{24}\omega_{35}\otimes
e_{5}e_{3}-r^{-1}s{\mathcal{E}}_{24}\omega_{35}\otimes e_{3}e_{5}) \\
&\quad + (r^{-1}e_{3}\omega_{245}\otimes
{\mathcal{E}}_{245}-r^{-1}e_{3}\omega_{245}\otimes
{\mathcal{E}}_{245}) +(\omega_{2453}\otimes
{\mathcal{E}}_{245}e_{3}-r^{-1}\omega_{2453}\otimes
e_{3}{\mathcal{E}}_{245})\\
&={\mathcal{E}}_{2453}\otimes 1 +
(1-r^{-1}s){\mathcal{E}}_{245}\omega_{3}\otimes e_{3} +
(1-r^{-1}s){\mathcal{E}}_{243}\omega_{5}\otimes e_{5}\\
&\quad +(1-r^{-1}s)^{2}{\mathcal{E}}_{24}\omega_{35}\otimes
e_{5}e_{3} +
(1-r^{-1}s)^2e_{2}\omega_{345}\otimes{\mathcal{E}}_{45}e_{3}\\
&\quad -(1-r^{-1}s)r^{-1} e_{2}\omega_{345}\otimes
{\mathcal{E}}_{345} +\omega_{2453}\otimes{\mathcal{E}}_{2453}.
\end{split}
\end{equation*}

\noindent{\bf Remark.} \ As indicated in Lemma 3.1,
 the right hand-side of the formula above does show that each product's ordering in the summation consisted of those possible good letters
(appearing as the 2nd factors in those tensor monomials) satisfies
the required non increasing property with respect to the ordering
$<$.

\vskip30pt \centerline{\bf ACKNOWLEDGMENT}

\vskip15pt NH would like to express his thanks to the ICTP for
support when he is visiting at the ICTP Mathematics Section,
Trieste, Italy, during March 1st to August 31st, 2006.

\end{document}